\documentclass[11pt]{article}

\usepackage[T1]{fontenc}
\usepackage[utf8]{inputenc}
\usepackage{lmodern}
\usepackage{amsmath,amssymb,amsthm,mathtools}
\usepackage[a4paper,margin=28mm]{geometry}
\usepackage{microtype}
\usepackage{enumitem}
\usepackage[hidelinks]{hyperref}
\usepackage{xcolor}
\usepackage{tikz}

\numberwithin{equation}{section}

\newtheorem{theorem}{Theorem}[section]
\newtheorem{lemma}[theorem]{Lemma}
\newtheorem{proposition}[theorem]{Proposition}
\newtheorem{corollary}[theorem]{Corollary}
\newtheorem{claim}{Claim}[theorem]

\theoremstyle{definition}
\newtheorem{definition}[theorem]{Definition}
\newtheorem{example}[theorem]{Example}
\newtheorem{remark}[theorem]{Remark}
\newtheorem{conjecture}[theorem]{Conjecture}

\newcommand{\rh}{\rho}
\newcommand{\spex}{\operatorname{spex}}
\newcommand{\SPEX}{\operatorname{SPEX}}
\DeclareMathOperator{\ex}{ex}
\DeclareMathOperator{\EX}{EX}
\newcommand{\Estar}{E_*}
\newcommand{\dst}{\delta_*}
\newcommand{\Pprof}{\mathcal P}
\newcommand{\FJ}{F_J}
\newcommand{\FL}{F_L}
\newcommand{\Jn}[2]{J_n(#1,#2)}
\newcommand{\Ln}[2]{L_n(#1,#2)}
\newcommand{\ST}[1]{\mathrm{ST}_2(#1)}
\newcommand{\Ttwo}[1]{T_2(#1)}
\newcommand{\Oh}{O_H}
\newcommand{\px}{\boldsymbol{x}}
\newcommand{\pxmax}{x_{\max}}
\newcommand{\pxmin}{x_{\min}}

\title{Non-bipartite spectral extremal graphs forbidding\\
a three-chromatic color-critical graph%
\thanks{Supported by the National Natural Science Foundation of China (Nos. 12271439 and 12301452), the China Scholarship Council (No. 202506290205), the Natural Science Foundation of Shanxi Province (No. 202303021212025), and Shanxi Scholarship Council of China (No. 2025-080).}}

\author{%
Bo~Dong$^{a,b,d}$,
Cunxiang~Duan$^{c}$,
Attila~Sali$^{d,e}$\thanks{Corresponding author.},
Ligong~Wang$^{a,b}$\\[1ex]
\parbox{0.94\textwidth}{%
\centering\small
$^{a}$School of Mathematics and Statistics,
Northwestern Polytechnical University,\\
Xi'an, Shaanxi 710129, P.R. China\\[0.5ex]
$^{b}$Xi'an-Budapest Joint Research Center for Combinatorics,\\
Northwestern Polytechnical University,
Xi'an, Shaanxi 710129, P.R. China\\[0.5ex]
$^{c}$School of Mathematics, Taiyuan University of Technology,\\
Taiyuan, Shanxi 030024, P.R. China\\[0.5ex]
$^{d}$Alfr\'ed R\'enyi Institute of Mathematics,
Budapest, Hungary\\[0.5ex]
$^{e}$Department of Computer Science and Information Theory,\\
Budapest University of Technology and Economics, Hungary\\[1ex]
E-mail: \texttt{dongbo98@163.com}; \texttt{cxduanmath@163.com};\\
\texttt{sali.attila@renyi.hu}; \texttt{lgwangmath@163.com}
}}

\date{}

\begin{document}

\maketitle

\begin{abstract}
A graph is color-critical if it contains at least one color-critical edge. For every fixed three-chromatic color-critical graph $H$, we determine the maximum spectral radius of an $n$-vertex non-bipartite $H$-free graph and characterize the unique extremal graph for all sufficiently large $n$.
The extremal graph belongs to one of three graph families:
$J_n(s,t)$, $\mathrm{ST}_2(n-1)$ or $L_n(s,t)$.
Explicit conditions on $H$ determine the family and its exact parameters.
This paper generalizes the results for books and theta graphs to all three-chromatic color-critical graphs.
Furthermore, an three-chromatic example demonstrates that spectral extremal graphs are not necessarily edge-extremal, disproving a conjecture of Fang and Lin [J. Algebraic Combin., 2026].
\end{abstract}

\section{Introduction}

All graphs are finite and simple.
A graph $G=(V(G),E(G))$ consists of the vertex set $V(G)$ and the edge set $E(G)$. The number of vertices (resp., edges) of $G$ is called its \emph{order} (resp., \emph{size}). We write $e(G)=|E(G)|$. For $v\in V(G)$, let $N_G(v)$ be its neighborhood and $d_G(v)=|N_G(v)|$ its degree.
Graph containment $H\subseteq G$ means that $G$ contains a subgraph isomorphic to $H$, and the subgraph need not be induced.
A graph $G$ is $H$-free if $H\nsubseteq G$.
In an embedding of $H$ into $G$, we call $G$ the host graph.

We write $A(G)$, $\rh(G)$, and $\chi(G)$ for the adjacency matrix, adjacency spectral radius, and chromatic number of a graph $G$, respectively. Throughout the paper, we abbreviate $\rh(G)$ to $\rho$ unless otherwise specified.
For a connected graph $G$, a positive eigenvector of $A(G)$ corresponding to $\rho(G)$ is called a \emph{Perron vector}; we always normalize it to have Euclidean norm one.

For a fixed graph $H$ and an integer $r\ge2$, define
\[
 \ex_{r+1}(n,H)
 =
 \max\{
 e(G): G \text{~is an $n$-vertex graph},\ H\nsubseteq G,\ \chi(G)>r
 \},
\]
and let $\EX_{r+1}(n,H)$ be the family of graphs attaining this maximum number of edges.
Similarly, define
\[
 \spex_{r+1}(n,H)
 =
 \max\{
 \rh(G): G\text{~is an $n$-vertex graph},\ H\nsubseteq G,\ \chi(G)>r
 \},
\]
and let $\SPEX_{r+1}(n,H)$ be the family of graphs attaining this maximum spectral radius.

An edge $u_1v_1$ of a graph $H$ is \emph{color-critical} if
$\chi(H-u_1v_1)<\chi(H)$.
For a three-chromatic graph, this is equivalent to $H-u_1v_1$ being bipartite.
Set
\[
 \Estar(H)
 =
 \{u_1v_1\in E(H):H-u_1v_1\text{ is bipartite}\},
 \quad
 \dst(H)
 =
 \min_{u_1v_1\in\Estar(H)}
 \min\{d_H(u_1),d_H(v_1)\}.
\]
A graph is called \emph{color-critical} if it has a color-critical edge.
Throughout the paper, fix a graph $H$ with $\chi(H)=3$
and $\Estar(H)\ne\varnothing$. Set $n_0=|V(H)|$.
The notation $O(\cdot)$ uses absolute implicit constants, whereas $O_H(\cdot)$ allows them to depend on $H$; neither depends on $n$. The $r$-partite Tur\'an graph $T_r(n)$ is the complete $r$-partite graph of order $n$ whose part sizes differ by at most one. The balanced complete bipartite graph of order $n$ is denoted by $\Ttwo n$, so $e(\Ttwo n)=\left\lfloor n^2/4\right\rfloor$.

The classical Tur\'an problem asks for the maximum number of edges in an $n$-vertex graph containing no prescribed subgraphs.
Tur\'an's theorem~\cite{Turan1941} identifies $T_r(n)$ as the unique extremal graph when the forbidden graph is $K_{r+1}$.
Simonovits' chromatic critical edge theorem~\cite{Simonovits1968} extends this conclusion, for sufficiently large $n$, to every fixed color-critical graph $H$ with $\chi(H)=r+1$.
Nikiforov established the spectral analogue of Tur\'an's theorem~\cite{Nikiforov2007}.
A spectral counterpart of Simonovits' critical-edge theorem follows from his work on spectral saturation~\cite{Nikiforov2009}, together with an analysis of the equality case; see Zhai and Lin~\cite[Theorem~1.2]{ZhaiLin2023} for an explicit statement including uniqueness.
Specifically, for every fixed color-critical graph $H$ with $\chi(H)=r+1$, the unique $n$-vertex $H$-free graph of maximum spectral radius is $T_r(n)$ for sufficiently large $n$.
Spectral extremal graph theory, especially spectral Tur\'an problems, has attracted considerable attention in recent years. This area has produced many significant results, for an overview of its development and major advances, see the surveys by Nikiforov~\cite{Nikiforov2011survey} and by Li, Liu, and Feng~\cite{LiLiuFeng2022survey}.

By Nikiforov's spectral color-critical edge theorem~\cite{Nikiforov2009} (see also Zhai and Lin~\cite[Theorem~1.2]{ZhaiLin2023}), for every fixed three-chromatic color-critical graph $H$ and all sufficiently large $n$, $\Ttwo n$ is the unique $n$-vertex $H$-free graph of maximum spectral radius.
Since $\Ttwo n$ is bipartite, this result naturally motivates the study of the corresponding problem for non-bipartite graphs: what is the maximum spectral radius among non-bipartite $H$-free graphs of order $n$, and which graphs attain it?
Lin, Ning, and Wu~\cite{LinNingWu2021} proved that, for $n\ge5$,
\[
 \spex_3(n,C_3)=\rho(\ST{n-1}),
 \qquad
 \SPEX_3(n,C_3)=\{\ST{n-1}\},
\]
where $\ST{n-1}$ is obtained by subdividing one edge of $T_2(n-1)$.
Li and Peng~\cite{LiPeng2023} subsequently gave an alternative proof and extended this refinement to non-$r$-partite $K_{r+1}$-free graphs.

For every fixed $k\ge2$ and sufficiently large $n$, Zhang and Zhao~\cite{ZhangZhao2023} proved that
\[
 \spex_3(n,C_{2k+1})=\rho(L_n(1,1)),
 \qquad
 \SPEX_3(n,C_{2k+1})=\{L_n(1,1)\},
\]
where $L_n(1,1)$ is obtained by identifying a vertex of a triangle with a vertex in a smallest part of $T_2(n-2)$.
More recently, Zou, Feng, and Li~\cite{ZouFengLi2026} obtained the explicit sufficient threshold $n\ge187k$ and proved a more general result for families of forbidden short odd cycles.

A book $B_{r+1}$ consists of $r+1$ triangles sharing one edge.
For books and certain theta graphs $H$, Zhai and
Lin~\cite{ZhaiLin2023} obtained linear order thresholds
ensuring that $T_2(n)$ is the unique $n$-vertex $H$-free
graph of maximum spectral radius, with both bipartite
and non-bipartite graphs allowed.
Liu and Miao~\cite{LiuMiao2025} proved that, for $r\ge1$ and $n\ge8(r^2+r+4)$,
\[
 \spex_3(n,B_{r+1})=\rho(J_n(r,r)),
 \qquad
 \SPEX_3(n,B_{r+1})=\{J_n(r,r)\},
\]
where $J_n(r,r)$ is obtained from $T_2(n-1)$ by adding one vertex with exactly $r$ neighbors in each part.

For a theta graph $\theta(1,q,r)$, the three internally vertex-disjoint paths joining its two common endpoints have lengths $1,q,r$.
Fang and Lin~\cite{FangLin2026} determined both the edge and spectral extremal graphs when $q$ is even and $q,r\ge2$.
In the notation defined below, their spectral conclusion is
\[
 \SPEX_3(n,\theta(1,q,r))=
 \begin{cases}
 \{J_n(1,1)\},&r\text{ is even},\\
 \{\ST{n-1}\},&r\text{ is odd and }q=2,\\
 \{L_n(1,1)\},&r\text{ is odd and }q\ge4,
 \end{cases}
\]
for sufficiently large $n$.
In every case, their spectral extremal graph is also edge-extremal. Therefore they gave following conjecture.

\begin{conjecture}\label{conj:edge-spectral}
Let $H$ be a fixed color-critical graph with $\chi(H)=r+1\ge3$.
Then, for all sufficiently large $n$,
\[
 \SPEX_{r+1}(n,H)\subseteq\EX_{r+1}(n,H).
\]
\end{conjecture}

However, Wu and Lu~\cite{WuLu2026} disproved this conjecture using
a forbidden graph of chromatic number five.
Our example shows that the same separation already occurs
for a connected triangle-free forbidden graph of
chromatic number three.

The preceding results raise a structural question: which properties of a three-chromatic color-critical graph determine its non-bipartite spectral extremal graph?
We prove that the minimum endpoint degree over the critical
edges of $H$, together with the presence or absence of a
triangle, determines which of the three graph families
defined below contains the spectral extremal graph.
The numbers of neighbors in the two parts after deleting a vertex of $H$, or the degrees of the endpoints of a critical edge, then we determine the exact parameters $s$ and $t$.
This gives a classification for every fixed three-chromatic color-critical graph, including graphs with leaves and disconnected graphs.
It extends the known book and theta results and provides a counterexample to the Fang--Lin conjecture.

\begin{definition}\label{def:families}
The following graphs have order $n$.

$(i)$ For $1\le t\le s\le\lfloor(n-1)/2\rfloor$, the graph
$\Jn st$ is obtained from $\Ttwo{n-1}$ by adding a vertex $u$ with $s$ neighbors in a smaller
part and $t$ neighbors in the other part.
These graphs form the $J$-family.

$(ii)$ The graph $\ST{n-1}$ is obtained by subdividing one edge of $\Ttwo{n-1}$.

$(iii)$ For $1\le t\le s\le\lfloor(n-2)/2\rfloor$, the graph
$\Ln st$ is obtained from $\Ttwo{n-2}$ by adding adjacent vertices $u,v$.
Choose sets $V_t^{(1)}\subseteq V_s^{(1)}$ in a smaller part, with
$|V_s^{(1)}|=s$ and $|V_t^{(1)}|=t$, and set $N(u)\setminus\{v\}=V_s^{(1)}$, $N(v)\setminus\{u\}=V_t^{(1)}$.
These graphs form the $L$-family.
\end{definition}

In the $J$ and $L$ constructions, the specified complete
bipartite graph is called the \emph{core}.
Its parts are denoted by $V^{(1)},V^{(2)}$, with $|V^{(1)}|\le|V^{(2)}|$.
When the two parts have equal size, either may be
designated as $V^{(1)}$.
The added vertex $u$ in $\Jn st$ and the added vertices
$u,v$ in $\Ln st$ are called \emph{exceptional vertices}.
Their neighbors in the core are called their
\emph{core neighbors}.
Thus, in $\Jn st$, the vertex $u$ has $s$ core neighbors
in $V^{(1)}$ and $t$ in $V^{(2)}$; in $\Ln st$, the vertices $u,v$
have nested core neighborhoods of sizes $s,t$,
respectively, both contained in $V^{(1)}$.


The next theorem gives the complete classification.
Subsections~\ref{sec:Jparameters} and~\ref{sec:Lparameters}
give explicit finite rules for determining the parameters.

\begin{theorem}\label{thm:main}
Let $H$ be a fixed three-chromatic color-critical graph.
For all sufficiently large $n$, the family
$\SPEX_3(n,H)$ consists of a single graph, as follows:
\[
 \SPEX_3(n,H)=
 \begin{cases}
 \{\Jn{s_J}{t_J}\},
 &\dst(H)\ge3,\\[2pt]
 \{\ST{n-1}\},
 &\dst(H)=2\text{ and }H\text{ contains a triangle},\\[2pt]
 \{\Ln{s_L}{t_L}\},
 &\dst(H)=2\text{ and }H\text{ is triangle-free}.
 \end{cases}
\]
The parameters $s_J\ge t_J\ge1$ are determined by
Theorem~\ref{j:exact-parameters}, and $s_L\ge t_L\ge1$
by Theorem~\ref{l:parameters}.
Consequently, $\spex_3(n,H)$ is the spectral radius
of the graph in the corresponding line.
\end{theorem}

The remainder of the paper is organized as follows.
Section~2 collects the preliminary results used in the
structural and spectral arguments.
Section~\ref{sec:structure} establishes the three possible extremal graph types: the $J$-family, the $L$-family, and $\ST{n-1}$.
Subsections~\ref{sec:Jparameters} and~\ref{sec:Lparameters} determine the exact parameters of the $J$ and $L$ families through containment criteria and spectral comparisons.
Finally, Section~\ref{sec:applications} reduces the problem to connected forbidden graphs and completes the proof of Theorem~\ref{thm:main}, applies the classification to color-critical generalized theta graphs, and presents asymmetric examples, including a counterexample to the Fang--Lin conjecture.

\section{Preliminaries}\label{sec:preliminaries}

In this section, we set up the basic notation and state the preliminary lemmas used in the sequel.

For a real vector $\boldsymbol y$, write $\|\boldsymbol y\| = \left(\sum_i y_i^2\right)^{1/2}.$
For a real matrix $M=(m_{ij})$, its Euclidean operator norm and Frobenius norm are, respectively,
\[
 \|M\|
 =
 \sup_{\|\boldsymbol y\|=1}\|M\boldsymbol y\|,
 \qquad
 \|M\|_{\mathrm F}
 =
 \left(\sum_{i,j}m_{ij}^2\right)^{1/2}.
\]
In particular, $\|M\|_{\mathrm F}^2=\operatorname{tr}(M^TM)$.
If $\lambda_1,\ldots,\lambda_n$ are the adjacency eigenvalues of a graph $G'$, the spectral theorem gives
\[
 \rh(G')^2
 \le\sum_{i=1}^n\lambda_i^2
 =\operatorname{tr}(A(G')^2)
 =\|A(G')\|_{\mathrm F}^2
 =2e(G').
\]

The following lemma establishes the elementary properties of critical edges required for every structural case.

\begin{lemma}\label{prelim:critical}
Let $H$ be three-chromatic and color-critical.
Then every critical edge of $H$ belongs to every odd cycle of $H$ and is not a bridge, and $\dst(H)\ge2$.
Suppose in addition that $u_1v_1\in\Estar(H)$ and $d_H(u_1)=2$, with $N_H(u_1)=\{v_1,v_2\}$.
Then $H-u_1$ is bipartite, $v_1$ and $v_2$ lie in opposite parts of the bipartition inherited from $H-u_1v_1$, and $H$ contains a triangle if and only if $v_1v_2\in E(H)$.
\end{lemma}

\begin{proof}
An odd cycle avoiding $u_1v_1$ would remain in $H-u_1v_1$.
Thus every odd cycle contains $u_1v_1$, and an odd cycle containing $u_1v_1$ shows that this edge is not a bridge.
In particular, both endpoints of every critical edge have degree at least two.

In a bipartition of $H-u_1v_1$, the vertices $u_1,v_1$ have the same color, since otherwise it would also give a proper two-coloring of $H$.
The vertex $v_2$ has the opposite color.
Every triangle must contain $u_1v_1$ and hence, when $d_H(u_1)=2$, it must be $u_1v_1v_2$.
\end{proof}


\begin{lemma}[\cite{OWu2026}]\label{prelim:reduction}
Let $H$ be a connected three-chromatic color-critical graph. For sufficiently large $n$, if $G\in\SPEX_3(n,H)$ and $\rh(G)>(n-3)/2$, then $G$ is connected.
Every vertex $u$ of minimum Perron coordinate satisfies $\chi(G-u)=2$.
\end{lemma}

This is the specialization $r=2$ of Remark~2.17 in~\cite{OWu2026}.

\begin{lemma}[\cite{LinNingWu2021}]\label{prelim:triangle}
For sufficiently large $n$, if $G$ is a non-bipartite triangle-free graph of order $n$, then
\[
 \rh(G)\le\rh(\ST{n-1}),
\]
with equality if and only if $G\cong\ST{n-1}$.
\end{lemma}

In Sections~2--3, assume in addition that $H$ is connected.
The extension to disconnected forbidden graphs will be established
in the proof of Theorem~\ref{thm:main}.
For a connected $H$-free graph $G$ of order $n$, let $\rho=\rh(G)$ and let $\px=(x_v)_{v\in V(G)}$ be its positive unit Perron vector and
\[
 \pxmax=\max_{v\in V(G)}x_v,
 \qquad
 \pxmin=\min_{v\in V(G)}x_v.
\]
Note that $\pxmin^2\le1/n$. Whenever $\rho\ge n/3$, the eigenvalue equation and Cauchy--Schwarz give
\[
 \rho\pxmax\le\sqrt n,
 \qquad
 \pxmax\le\frac3{\sqrt n}.
\]

The next lemma measures the spectral cost of deleting vertices
with small Perron coordinates.

\begin{lemma}\label{prelim:delete}\label{prelim:loss}
Let $G'$ be a connected graph of order $n\ge2$, let $\rho=\rh(G')$,
and let $\px$ be its positive unit Perron vector.
For $X\subseteq V(G')$, let $\vartheta=\sum_{v\in X}x_v^2$.

$(i)$ If $\vartheta<1$, then
\[
 \rh(G'-X)\ge
 \frac{
 \rho(1-2\vartheta)
 +2\displaystyle\sum_{vw\in E(G'[X])}x_vx_w
 }{1-\vartheta}.
\]

$(ii)$ If $\vartheta\le1/4$, then
\[
 \rho^2-\rh(G'-X)^2\le3n^2\vartheta.
\]

$(iii)$ If $u$ is a minimum-coordinate vertex, then
\[
 \rh(G'-u)\ge\rho\frac{n-2}{n-1},
 \qquad
 \rho^2-\rh(G'-u)^2<2n.
\]
\end{lemma}

\begin{proof}
$(i).$ Multiplying the eigenvalue equation at each $v\in X$
by $x_v$ and summing over $X$ yields
\[
 \rho\vartheta
 =
 2\sum_{vw\in E(G'[X])}x_vx_w
 +
 \sum_{\substack{v\in X,\ w\in V(G')\setminus X\\
                  vw\in E(G')}}x_vx_w.
\]
Consequently, the adjacency quadratic form remaining
on $V(G')\setminus X$ equals
\[
 \rho(1-2\vartheta)
 +2\sum_{vw\in E(G'[X])}x_vx_w.
\]
The restricted vector $\left.\px\right|_{V(G'-X)}$
has squared norm $1-\vartheta$, so the first inequality
follows from the Rayleigh principle.

$(ii).$ Discarding the nonnegative internal-edge term gives $\rh(G'-X)\ge
 \rho\frac{1-2\vartheta}{1-\vartheta}.$
For $\vartheta\le1/4$, the right-hand side is nonnegative and
\[
 1-
 \left(\frac{1-2\vartheta}{1-\vartheta}\right)^2
 =
 \frac{2\vartheta-3\vartheta^2}{(1-\vartheta)^2}
 \le3\vartheta.
\]
Since $\rho\le n-1\le n$, this proves $\rho^2-\rh(G'-X)^2\le3n^2\vartheta.$

$(iii)$. For a minimum-coordinate vertex $u$, we have
$\vartheta=x_u^2\le1/n$.
The function $(1-2\vartheta)/(1-\vartheta)$ is decreasing
for $\vartheta<1$, so the first inequality with $X=\{u\}$
gives $\rh(G'-u)\ge\rho\frac{n-2}{n-1}.$
The squared loss is therefore at most $\frac{\rho^2(2n-3)}{(n-1)^2}
 \le2n-3<2n.$
\end{proof}

\begin{lemma}\label{prelim:bipartite}
If $B$ is a bipartite graph with part sizes $p,q$, then
\[
 \rh(B)^2\le e(B)\le pq.
\]
\end{lemma}

\begin{proof}
Let $A_{12}$ be the $p\times q$ cross-adjacency matrix of $B$.
Then
\[
 A(B)=
 \begin{pmatrix}
  0&A_{12}\\
  A_{12}^T&0
 \end{pmatrix},
 \qquad
 A(B)^2=
 \begin{pmatrix}
  A_{12}A_{12}^T&0\\
  0&A_{12}^TA_{12}
 \end{pmatrix}.
\]
Since $A_{12}A_{12}^T$ is positive semidefinite,
\[
 \rh(B)^2
 =
 \rh(A_{12}A_{12}^T)
 \le\operatorname{tr}(A_{12}A_{12}^T)
 =
 \|A_{12}\|_{\mathrm F}^2
 =
 e(B).
\]
The last equality holds because $A_{12}$ records each cross-edge once.
There are at most $pq$ such edges, so $e(B)\le pq$.
\end{proof}

The Schur-complement, resolvent, Rayleigh principle and Perron--Frobenius results are standard; see~\cite{BhatiaDavis,BrouwerHaemers,HornJohnson}.

\begin{lemma}[\cite{BhatiaDavis}]\label{prelim:residual}
For a graph $G'$ of order $n$, let $\boldsymbol y$ be a unit vector and $r=\boldsymbol y^TA(G')\boldsymbol y$.
Then
\[
 \rh(G')-r
 \ge
 \frac{\|A(G')\boldsymbol y-r\boldsymbol y\|^2}{2n}.
\]
\end{lemma}

Indeed, the Bhatia--Davis variance bound~\cite{BhatiaDavis} gives, for a symmetric matrix with largest and smallest eigenvalues $\lambda_1,\lambda_n$, $\|A\boldsymbol y-r\boldsymbol y\|^2
 \le(\lambda_1-r)(r-\lambda_n).$
For every $n$-vertex graph, $r-\lambda_n\le2n$, which yields the displayed consequence.

\begin{definition}\label{def:schur-complement}
Let
\[
 A=
 \begin{pmatrix}
  A_{11}&A_{12}\\
  A_{21}&A_{22}
 \end{pmatrix}
\]
be a square block matrix with $A_{11}$ invertible.
The \emph{Schur complement of $A_{11}$ in $A$} is
\[
 A/A_{11}
 :=
 A_{22}-A_{21}A_{11}^{-1}A_{12}.
\]
The same notation is used for any invertible principal diagonal block
after a simultaneous permutation of the corresponding block rows and columns.
\end{definition}

\begin{lemma}[\cite{HornJohnson,BrouwerHaemers}]
\label{lem:schur-determinant}
Let $A$ satisfy the hypotheses of Definition~\ref{def:schur-complement}.
Then
\[
 \det(A)=\det(A_{11})\det(A/A_{11}).
\]
\end{lemma}

\begin{definition}\label{def:core-effective}
For a graph $G'$, let $X\subseteq V(G')$ be the vertex set
whose coordinates are retained when those indexed by
$V(G')\setminus X$ are eliminated.
With the vertices of $X$ listed first, write
\[
 A(G')=
 \begin{pmatrix}
  A_{11}&A_{12}\\
  A_{12}^T&A_{22}
 \end{pmatrix}.
\]
Here $A_{11}$ is the adjacency matrix of $G'[X]$,
$A_{22}$ is the adjacency matrix of $G'-X$, and $A_{12}$
records the adjacencies between $X$ and $V(G')\setminus X$.
For $\rho>\rh(A_{22})$, define
\[
 R_X(\rho):=(\rho I-A_{22})^{-1},
 \qquad
 K_X(\rho):=A_{11}+A_{12}R_X(\rho)A_{12}^T.
\]
We call $R_X(\rho)$ the \emph{resolvent associated with $X$} and $K_X(\rho)$ the \emph{effective matrix on $X$}.
\end{definition}

The next lemma expresses the eigenvalue equation solely
in the coordinates indexed by $X$, after eliminating
those indexed by $V(G')\setminus X$.

\begin{lemma}[\cite{BrouwerHaemers,HornJohnson}]
\label{lem:core-elimination}
With the notation of Definition~\ref{def:core-effective},
let $\rho>\rh(A_{22})$.
Then
\[
 A(G')
 \begin{pmatrix}
  \boldsymbol z\\
  \boldsymbol y
 \end{pmatrix}
 =
 \rho
 \begin{pmatrix}
  \boldsymbol z\\
  \boldsymbol y
 \end{pmatrix}
\]
is equivalent to
\[
 \boldsymbol y=R_X(\rho)A_{12}^T\boldsymbol z,
 \qquad
 K_X(\rho)\boldsymbol z=\rho\boldsymbol z.
\]
Moreover,
\[
 \det(\rho I-A(G'))
 =
 \det(\rho I-A_{22})
 \det(\rho I-K_X(\rho)).
\]
If $G'$ is connected, $\rho=\rh(G')$, and
$\begin{pmatrix}\boldsymbol z\\\boldsymbol y\end{pmatrix}$
is a positive Perron vector of $G'$, then $\boldsymbol z$
is positive, $K_X(\rho)$ is nonnegative, and
$\rho=\rh(K_X(\rho))$.
\end{lemma}

The determinant identity in Lemma~\ref{lem:core-elimination} follows
from Lemma~\ref{lem:schur-determinant} applied to $\rho I-A(G')$
after interchanging the two block rows and the corresponding
block columns.

\begin{lemma}[\cite{HornJohnson}]
\label{lem:resolvent-neumann}
With the notation of Definition~\ref{def:core-effective},
let $\rho>\rh(A_{22})$.
Then
\[
 R_X(\rho)
 =
 \frac1\rho
 \sum_{j=0}^{\infty}
 \left(\frac{A_{22}}{\rho}\right)^j.
\]
In particular, $R_X(\rho)$ is entrywise nonnegative.

More generally, if a square matrix $M$ satisfies $\|M\|<1$, then
\[
 (I-M)^{-1}
 =
 I+M+M^2(I-M)^{-1},
\]
and
\[
 \|(I-M)^{-1}\|
 \le\frac1{1-\|M\|},
 \qquad
 \|(I-M)^{-1}-I-M\|
 \le\frac{\|M\|^2}{1-\|M\|}.
\]
\end{lemma}

The next lemma supplies the lower bounds needed before invoking the deletion principle.

\begin{lemma}\label{prelim:lower}
Let $H$ be a fixed connected three-chromatic color-critical graph.
If $n$ is sufficiently large, then the following hold.
\begin{enumerate}[label=(\roman*),leftmargin=*]
\item
If $\dst(H)\ge3$, then $\Jn11$ is $H$-free and
$
 \rh(\Jn11)^2>e(\Ttwo{n-1}).
$

\item
If $\dst(H)=2$ and $H$ contains a triangle, then $\ST{n-1}$ is $H$-free and
$
 \rh(\ST{n-1})^2>e(\Ttwo{n-1})-2.
$

\item
If $\dst(H)=2$ and $H$ is triangle-free, then $\Ln11$ is $H$-free and
$
 \rh(\Ln11)^2>e(\Ttwo{n-2}).
$
\end{enumerate}
In every case,
$
 \spex_3(n,H)>\frac{n-3}{2},~
 \spex_3(n,H)>\frac n3,~
 \spex_3(n,H)^2>e(\Ttwo{n-2}).
$
\end{lemma}

\begin{proof}
$(i).$ Any copy of $H$ in $\Jn11$ must use the added degree-two vertex, since deleting this vertex leaves a bipartite graph.
Its preimage $v_1\in V(H)$ therefore satisfies $H-v_1$ bipartite and $d_H(v_1)=2$.
The two neighbors of $v_1$ lie in opposite core parts, so either edge incident with $v_1$ is critical.
This contradicts $\dst(H)\ge3$.
Thus $\Jn11$ is $H$-free.
Since $\Ttwo{n-1}$ is a proper induced subgraph of the connected graph $\Jn11$,
\[
 \rh(\Jn11)^2
 >
 \rh(\Ttwo{n-1})^2
 =
 e(\Ttwo{n-1}).
\]

$(ii).$ $\ST{n-1}$ is $H$-free because it is triangle-free.
Let $(V^{(1)},V^{(2)})$ be the bipartition of $\Ttwo{n-1}$, with sizes $p=\left\lfloor\frac{n-1}{2}\right\rfloor,$ $q=\left\lceil\frac{n-1}{2}\right\rceil.$
Fix an edge $v_1^{(1)}v_1^{(2)}$, where $v_1^{(i)}\in V^{(i)}$ for $i=1,2$.
By symmetry, its unit Perron vector $\boldsymbol x$ satisfies
$2x_{v_1^{(1)}}x_{v_1^{(2)}}=1/\sqrt{pq}$.
Evaluating the Rayleigh quotient at $\boldsymbol x$ gives
\[
 \rh(\Ttwo{n-1}-v_1^{(1)}v_1^{(2)})
 \ge
 \sqrt{pq}-\frac1{\sqrt{pq}}.
\]
Since $\Ttwo{n-1}-v_1^{(1)}v_1^{(2)}$ is a proper subgraph of the connected graph $\ST{n-1}$,
\[
 \rh(\ST{n-1})^2
 >
 \rh(\Ttwo{n-1}-v_1^{(1)}v_1^{(2)})^2
 >
 pq-2
 =
 e(\Ttwo{n-1})-2.
\]

$(iii).$ $\Ln11$ is $H$-free because its only odd cycle is a triangle, whereas every odd cycle of $H$ has length at least five.
Moreover,
\[
 \rh(\Ln11)^2>e(\Ttwo{n-2})
\]
because $\Ttwo{n-2}$ is a proper induced subgraph of the connected graph $\Ln11$.

The remaining inequalities follow from these bounds and $e(\Ttwo m)=\lfloor m^2/4\rfloor$.
\end{proof}

In all three cases, fix $G\in\SPEX_3(n,H)$, and choose a vertex $u$ of minimum Perron coordinate. Lemmas~\ref{prelim:reduction} and~\ref{prelim:lower} show that $G$ is connected and $B=G-u$ is bipartite. Let $(V^{(1)},V^{(2)})$ be its bipartition, with $|V^{(1)}|=p$ and $|V^{(2)}|=q$.

For a bipartite graph $B$ with specified parts $V^{(1)},V^{(2)}$, our barred notation is relative to the complete bipartite host:
\[
 \bar E(B)
 =
 \{v^{(1)}v^{(2)}:v^{(1)}\in V^{(1)},\ v^{(2)}\in V^{(2)},\
 v^{(1)}v^{(2)}\notin E(B)\},
 \qquad
 \bar e(B)=|V^{(1)}||V^{(2)}|-e(B).
\]
For $v^{(1)}\in V^{(1)}$, write
\[
 \bar N_B(v^{(1)})=V^{(2)}\setminus N_B(v^{(1)}),
 \qquad
 \bar d_B(v^{(1)})=|\bar N_B(v^{(1)})|,
\]
and define the analogous notation for $v^{(2)}\in V^{(2)}$.
Thus $\bar d_B$ counts missing cross-edges.
For a graph without a specified bipartite host, barred notation has the usual meaning relative to the complete graph.

The following rough estimates are sufficient for the first structural arguments.

\begin{lemma}\label{prelim:rough}
Let $G$, $u$, $B$, and $(V^{(1)},V^{(2)})$ be as above. If $n$ is sufficiently large, then
\[
 e(B)>e(\Ttwo{n-2})-2n,
 \qquad
 \bar e(B)<3n,
 \qquad
 p,q=\frac n2+O(\sqrt n).
\]
Both $N_G(u)\cap V^{(1)}$ and $N_G(u)\cap V^{(2)}$ are nonempty.
\end{lemma}

\begin{proof}
By Lemmas~\ref{prelim:delete}, \ref{prelim:bipartite}, and~\ref{prelim:lower}, we have
\[
 e(B)\ge\rh(B)^2
 >
 \rho^2-2n
 >
 e(\Ttwo{n-2})-2n.
\]
Since $p+q=n-1$ and $pq\le e(\Ttwo{n-1})$, this gives
\[
 \bar e(B)
 <
 pq-\bigl(e(\Ttwo{n-2})-2n\bigr)
 <
 3n.
\]
The identity
\[
 pq=\frac{(n-1)^2-(p-q)^2}{4}
\]
implies $|p-q|=O(\sqrt n)$.
If $u$ had neighbors in only one part, then $G$ would be bipartite, a contradiction.
\end{proof}

The following lemma gives an upper bound on the degree of at least one of two neighbors of a prospective image of a degree-two critical endpoint.

\begin{lemma}\label{prelim:rectangle}
For a fixed three-chromatic color-critical graph $H$, assume that
$n_0=|V(H)|$, $u_1v_1\in\Estar(H)$, $d_H(u_1)=2$, and
$N_H(u_1)=\{v_1,v_2\}$.
Let $G'$ be an $H$-free graph such that $B=G'-u$ is bipartite
with parts $V^{(1)},V^{(2)}$.
Choose
\[
 v_1^{(1)}\in N_{G'}(u)\cap V^{(1)},
 \qquad
 v_1^{(2)}\in N_{G'}(u)\cap V^{(2)}.
\]
If $v_1v_2\in E(H)$, assume additionally that
$v_1^{(1)}v_1^{(2)}\in E(B)$.
Then
\[
 \min\{d_B(v_1^{(1)}),d_B(v_1^{(2)})\}
 \le n_0\sqrt{\bar e(B)+1}+1.
\]
\end{lemma}

\begin{proof}
Suppose that both degrees exceed the asserted bound.
Let
\[
 P=N_B(v_1^{(1)})\setminus\{v_1^{(2)}\},
 \qquad
 Q=N_B(v_1^{(2)})\setminus\{v_1^{(1)}\}.
\]
Then $|P||Q|>\bar e(B)n_0^2$.
Choose independently and uniformly random $n_0$-element subsets $P'\subseteq P$ and $Q'\subseteq Q$.
Each missing cross-edge between $P$ and $Q$ has both endpoints selected with probability $n_0^2/(|P||Q|)$.
If $Z$ counts the missing edges between $P'$ and $Q'$, then
\[
 \mathbb E[Z]
 \le
 \bar e(B)\frac{n_0^2}{|P||Q|}
 <
 1.
\]
Thus some choice gives a complete bipartite graph between $P'$ and $Q'$.

By Lemma~\ref{prelim:critical}, $H-u_1$ is bipartite with $v_1,v_2$ in opposite parts.
Map $u_1,v_1,v_2$ to $u,v_1^{(1)},v_1^{(2)}$, respectively.
Map the remaining vertices in the part containing $v_1$ into $Q'$, and those in the part containing $v_2$ into $P'$.
All required edges are present, including $v_1^{(1)}v_1^{(2)}$ when $v_1v_2\in E(H)$.
This embeds $H$ in $G'$, a contradiction.
\end{proof}

\section{Structural classification and exact parameters}
\label{sec:structure}\label{sec:parameters}

In this section, we characterize the structure and the exact parameters of the spectral extremal graphs for connected $H$. Recall that $n_0=|V(H)|$.


\subsection{The case $\dst(H)\ge3$: the $J$-family}
\label{sec:Jparameters}

In this subsection,  we show that the extremal graph belongs to the \(J\)-family and determine it exactly.

The following lemma bounds the degree of a minimum-coordinate vertex.


\begin{lemma}\label{j:bounded-degree}
If $\dst(H)\ge3$, then $d_G(u)=\Oh(1)$.
\end{lemma}

\begin{proof}
We first establish two local claims needed in the residual argument.
Let $\eta=1/(100n_0)$.
A core vertex $v\in V(B)$ is \emph{normal} if $\bar d_B(v)\le\eta n$, and \emph{abnormal} otherwise.
Write $V_{\mathrm{norm}}$ and $V_{\mathrm{abn}}$ for these two sets, and write $V_{\mathrm{norm}}^{(i)}=V_{\mathrm{norm}}\cap V^{(i)}$ and $V_{\mathrm{abn}}^{(i)}=V_{\mathrm{abn}}\cap V^{(i)}$ for $i\in\{1,2\}$.
By Lemma~\ref{prelim:rough}, $\bar e(B)<3n$, so
\[
 6n>
 \sum_{v\in V(B)}\bar d_B(v)
 \ge
 \sum_{v\in V_{\mathrm{abn}}}\bar d_B(v)
 \ge
 \frac{|V_{\mathrm{abn}}|n}{100n_0}.
\]
Consequently, $B$ has at most $K=600n_0$ abnormal vertices.

\begin{claim}\label{j:normal-coordinates}
If $n$ is sufficiently large and $v\in V_{\mathrm{norm}}$, then
$x_v\ge\pxmax/2$.
\end{claim}

\begin{proof}[Proof of the claim]
\renewcommand{\qedsymbol}{\openbox}
We first show that no abnormal core vertex attains the maximum Perron coordinate.
By Lemmas~\ref{prelim:lower} and~\ref{prelim:rough}, an abnormal core vertex $v$ satisfies
\[
 d_G(v)
 =
 \frac n2+O(\sqrt n)-\bar d_B(v)
 \le
 \frac n2-\eta n+O(\sqrt n)
 <
 \rho.
\]
However, the eigenvalue equation at a maximum-coordinate vertex $v$ implies $\rho\le d_G(v)$.

A maximum-coordinate vertex can always be chosen in $B$: if $u$ attains the maximum, its minimality forces all coordinates to be equal, so every vertex of $B$ also attains the maximum.
Thus a normal core vertex attains $\pxmax$.

For normal vertices $u_1^{(i)},u_2^{(i)}\in V^{(i)}$, where $i\in\{1,2\}$, their neighborhoods differ in at most $2\eta n+1$ vertices.
Hence
\[
 \rho|x_{u_1^{(i)}}-x_{u_2^{(i)}}|
 \le
 (2\eta n+1)\pxmax.
\]
Suppose that a maximum-coordinate normal vertex lies in $V^{(1)}$, and choose a normal vertex $v_1^{(2)}\in V^{(2)}$ with the maximum coordinate among the normal vertices of $V^{(2)}$.
The eigenvalue equation at the maximum-coordinate vertex in $V^{(1)}$ gives
\[
 \rho\pxmax
 \le
 |V^{(2)}|x_{v_1^{(2)}}+(K+1)\pxmax.
\]
Therefore
\[
 x_{v_1^{(2)}}
 \ge
 \frac{\rho-K-1}{|V^{(2)}|}\pxmax
 \ge
 \frac{n/2-K-2}{|V^{(2)}|}\pxmax
 =
 \bigl(1-\Oh(n^{-1/2})\bigr)\pxmax.
\]

The same-part estimate now applies in both parts.
For a normal vertex $u^{(1)}\in V^{(1)}$,
\[
 x_{u^{(1)}}
 \ge
 \left(1-\frac{2\eta n+1}{\rho}\right)\pxmax,
\]
while for a normal vertex $u^{(2)}\in V^{(2)}$,
\[
 x_{u^{(2)}}
 \ge
 \left(
 1-\frac{2\eta n+1}{\rho}-\Oh(n^{-1/2})
 \right)\pxmax.
\]
Since $\rho>n/2-1$ and $\eta=1/(100n_0)$, all normal coordinates exceed $\pxmax/2$ for large $n$.
\end{proof}

\begin{claim}\label{j:common-normal-neighbors}
If $v^{(2)}\in N_G(u)\cap V^{(2)}$, then
$|N_G(u)\cap N_B(v^{(2)})\cap V_{\mathrm{norm}}^{(1)}|
 \le n_0-1.$
The corresponding statement holds with $V^{(1)},V^{(2)}$ interchanged.
\end{claim}

\begin{proof}[Proof of the claim]
\renewcommand{\qedsymbol}{\openbox}
Suppose that $N_G(u)\cap N_B(v^{(2)})\cap V_{\mathrm{norm}}^{(1)}$ contains a set $S$ of $n_0$ normal vertices.
Their common neighborhood among normal vertices of $V^{(2)}$ has size at least
\[
 |V^{(2)}|-K-n_0\eta n>n_0+1
\]
for sufficiently large $n$, using the part-size estimate in Lemma~\ref{prelim:rough}.

Choose a critical edge $u_1v_1$ of $H$.
In a bipartition of $H-u_1v_1$, the vertices $u_1,v_1$ have the same color.
Map them to $u,v^{(2)}$, respectively.
Map all vertices of the opposite color into $S$, and the remaining vertices into the common neighborhood of $S$ in $V^{(2)}$, avoiding $v^{(2)}$.
This embeds $H$ in $G$, a contradiction.
\end{proof}

We now apply the residual estimate to bound $d_G(u)$.
Choose $u_1^{(1)}\in N_G(u)\cap V^{(1)}, v_1^{(2)}\in N_G(u)\cap V^{(2)}.$
Let $G'$ be obtained by completing the bipartite core $B$ and retaining only $uu_1^{(1)},uv_1^{(2)}$ at $u$.
It contains the triangle $uu_1^{(1)}v_1^{(2)}$ and is $H$-free by the degree-two-vertex argument in the proof of Lemma~\ref{prelim:lower}(i).

Let $d=d_G(u), C_H=K+2n_0, \Delta=\sum_{v^{(1)}v^{(2)}\in\bar E(B)}x_{v^{(1)}}x_{v^{(2)}}.$
Define $S^{(1)}=(N_G(u)\cap V_{\mathrm{norm}}^{(1)})\setminus(\{u_1^{(1)}\}\cup N_B(v_1^{(2)})),$
and $S^{(2)}=(N_G(u)\cap V_{\mathrm{norm}}^{(2)})\setminus(\{v_1^{(2)}\}\cup N_B(u_1^{(1)})).$
By Claim~\ref{j:common-normal-neighbors}, at most $n_0$ normal neighbors in each part are excluded from these sets.
There are at most $K$ abnormal core vertices.
Thus at most $C_H$ neighbors are omitted from $S^{(1)}\cup S^{(2)}$.

The missing edges $v^{(1)}v_1^{(2)}$ with $v^{(1)}\in S^{(1)}$ and $u_1^{(1)}v^{(2)}$ with $v^{(2)}\in S^{(2)}$ form disjoint collections.
Minimality of $x_u$ gives
\begin{align*}
 \Delta
 \ge
 \sum_{v^{(1)}\in S^{(1)}}x_{v^{(1)}}x_{v_1^{(2)}}
 +
 \sum_{v^{(2)}\in S^{(2)}}x_{u_1^{(1)}}x_{v^{(2)}}\ge
 x_u\sum_{v\in S^{(1)}\cup S^{(2)}}x_v\ge
 x_u
 \left(
 \sum_{v\in N_G(u)}x_v-C_H\pxmax
 \right).
\end{align*}

Let $r=\px^TA(G')\px$.
Then
\begin{align*}
 r
 &=
 \rho
 -
 2x_u
 \sum_{v\in N_G(u)\setminus\{u_1^{(1)},v_1^{(2)}\}}x_v
 +
 2\Delta\ge
 \rho-2C_Hx_u\pxmax\ge
 \rho-\frac{6C_Hd\pxmax^2}{n}.
\end{align*}
Since $\rho x_u\le d\pxmax$, $\rho\ge n/3$, $d\le n-1$ and $\pxmax^2\le9/n$, we have $r\ge\rho/2$ for large $n$.

At least $d-K$ neighbors of $u$ are normal.
Claim~\ref{j:normal-coordinates} gives
$\rho x_u
 \ge
 \frac{d-K}{2}\pxmax, rx_u
 \ge
 \frac{d-K}{4}\pxmax.$
On the other hand,
$(A(G')\px)_u=x_{u_1^{(1)}}+x_{v_1^{(2)}}\le2\pxmax.$
If $d\ge2K+16$, then
\[
 |(A(G')\px)_u-rx_u|
 \ge
 \frac{d-K-8}{4}\pxmax
 \ge
 \frac{d\pxmax}{8}.
\]
Lemma~\ref{prelim:residual} yields
$\rh(G')-r
 \ge
 \frac{d^2\pxmax^2}{128n}.$
Consequently,
$\rh(G')-\rho
 \ge
 \frac{\pxmax^2}{n}
 \left(\frac{d^2}{128}-6C_Hd\right).$
For $d>768C_H$, the right-hand side is positive, contradicting spectral extremality.
Hence $d_G(u)=\Oh(1)$.
\end{proof}

The next theorem turns the bounded-degree conclusion into the exact core structure required for the $J$-family.

\begin{theorem}\label{j:complete-core}
Let $H$ be fixed and connected with $\chi(H)=3$ and $\dst(H)\ge3$, and let $G\in\SPEX_3(n,H)$.
If $u$ has minimum Perron coordinate, then, for sufficiently large $n$, $G-u=\Ttwo{n-1}$.
Moreover, both core-neighbor counts $|N_G(u)\cap V^{(1)}|$ and $|N_G(u)\cap V^{(2)}|$ are positive and at most $\dst(H)-2$.
\end{theorem}

\begin{proof}
By Lemma~\ref{j:bounded-degree}, we have $d_G(u)=\Oh(1)$ and
$x_u
 \le
 \frac{d_G(u)\pxmax}{\rho}
 =
 \Oh(n^{-3/2}).$
Lemmas~\ref{prelim:delete} and~\ref{prelim:lower} give
$\rh(G-u)^2
 >
 e(\Ttwo{n-1})-\Oh(n^{-1}).$
By Lemma~\ref{prelim:bipartite}, we have
$\rh(G-u)^2
 \le
 e(G-u)
 \le
 e(\Ttwo{n-1}).$
The edge count is an integer, so for large $n$ it equals $e(\Ttwo{n-1})$.
Therefore the core is complete and balanced.

Choose a critical edge $u_1v_1$ with $d_H(u_1)=\dst(H)$.
A bipartition of $H-u_1v_1$ places $u_1,v_1$ in the same part.
After deleting $u_1$, its neighborhood meets one inherited bipartition class in one vertex and the other in $\dst(H)-1$ vertices.
If one core-neighbor count were at least $\dst(H)-1$, the positivity of the other core-neighbor count and the complete large core would allow an embedding with $u_1$ mapped to $u$.
This contradicts $H$-freeness.
Thus both core-neighbor counts lie between $1$ and $\dst(H)-2$.
\end{proof}

We now determine the preferred orientation and the unique pair of core-neighbor counts. Let $m=\dst(H)-2$.
For the $J$-family, define the set of neighbor-count pairs $\Pprof_J(H)$ as follows.
For each vertex $u_1$ such that $H-u_1$ is bipartite, and each ordered bipartition $(V^{(1)},V^{(2)})$ of $H-u_1$, include
\[
 \bigl(|N_H(u_1)\cap V^{(1)}|,\ |N_H(u_1)\cap V^{(2)}|\bigr)
\]
in $\Pprof_J(H)$.
All componentwise choices of bipartition are allowed.
This definition applies to every three-chromatic color-critical graph, including disconnected graphs, although the extremal arguments in this subsection assume that $H$ is connected.
Thus $\Pprof_J(H)$ is finite and symmetric under interchanging its coordinates.

We call these pairs the \emph{profiles} of $H$ for the $J$-family.
This term refers only to the two neighbor counts in the displayed definition. For positive integers $t$, define
\begin{equation*}
 \FJ(t)
 =
 \min\{\alpha-1:
       (\alpha,\beta)\in\Pprof_J(H),\ \beta\le t\}.
\end{equation*}
The set over which the minimum is taken is nonempty: an endpoint of a critical edge with degree $\dst(H)$ gives the profile $(\dst(H)-1,1)$.

The following lemma determines exactly which pairs of core-neighbor counts avoid $H$.

\begin{lemma}\label{j:frontier}
Let $s\ge t\ge1$ be fixed integers.
Then, for sufficiently large $n$, $H\subseteq\Jn st$ if and only if there exists $(\alpha,\beta)\in\Pprof_J(H)$ such that $\alpha\le s$ and $\beta\le t$.
Equivalently, $\Jn st$ is $H$-free if and only if $s\le\FJ(t)$.

The function $\FJ$ is nonincreasing and $\FJ(1)=m$.
Every admissible pair with $s\ge t\ge1$ satisfies $s,t\le m$.
For each feasible $t$, the largest spectral radius occurs at $s=\FJ(t)$.
\end{lemma}

\begin{proof}
Every copy of $H$ in $\Jn st$ must use the exceptional vertex $u$, because its deletion leaves a bipartite graph.
If its preimage is $u_1\in V(H)$, the core bipartition induces a profile $(\alpha,\beta)\in\Pprof_J(H)$ satisfying $\alpha\le s$ and $\beta\le t$.
Conversely, given such a profile, map its deleted vertex to $u$ and its neighbors into the two core neighborhoods of sizes $s$ and $t$.
The remaining vertices then embed in the sufficiently large complete bipartite core.

Monotonicity follows directly from the minimum.
Both coordinates of every profile are positive, since otherwise $H$ would be bipartite.
If a profile has $\beta=1$, the edge joining the deleted vertex to its unique neighbor in the second color class is critical.
The deleted vertex consequently has degree at least $\dst(H)$, so $\alpha\ge\dst(H)-1$.
A critical endpoint of degree $\dst(H)$ attains this bound, proving $\FJ(1)=m$.

Thus $s\le\FJ(t)\le m$, and $t\le s\le m$.
Increasing $s$ while keeping $t$ fixed adds edges to a connected graph, strictly increasing its spectral radius.
\end{proof}

Define the finite index set
$I_J=\{t:1\le t\le m,\ t\le\FJ(t)\}.$
It is nonempty because $1\in I_J$.
Set
\begin{equation*}
 T_J=
 \left\{t\in I_J:
 \FJ(t)+t=\max_{i\in I_J}\bigl(\FJ(i)+i\bigr)
 \right\}.
\end{equation*}

The following theorem identifies the unique optimal boundary pair: we first maximize $s+t$, then use a parity-dependent rule to resolve ties.

\begin{theorem}\label{j:exact-parameters}
Suppose that $\dst(H)\ge3$.
For sufficiently large $n$, define
\[
 t_J=
 \begin{cases}
 \max T_J,&n\text{~is odd},\\
 \min T_J,&n\text{~is even},
 \end{cases}
 \qquad
 s_J=\FJ(t_J).
\]
Then
\[
 \SPEX_3(n,H)=\{\Jn{s_J}{t_J}\}.
\]
\end{theorem}

\begin{proof}
We first derive the scalar spectral equation and its uniform expansion.
Let
\[
 p=\left\lfloor\frac{n-1}{2}\right\rfloor,
 \qquad
 q=\left\lceil\frac{n-1}{2}\right\rceil.
\]

\begin{claim}\label{j:spectral-equation}
For fixed integers $s\ge t\ge1$, let $\sigma=s+t$.
Then $\rh(\Jn st)$ is the unique solution greater than
$\sqrt{pq}$ of
\begin{equation}\label{j:spectral-resolvent}
 \rho
 =
 \frac{\sigma}{\rho}
 +
 \frac{qs^2+pt^2+2\rho st}
      {\rho(\rho^2-pq)}.
\end{equation}
Equivalently,
\begin{equation}\label{j:spectral-polynomial}
 (\rho^2-\sigma)(\rho^2-pq)
 =
 qs^2+pt^2+2\rho st.
\end{equation}

For fixed core-neighbor counts $s\ge t$, assigning $s$ to the smaller core part maximizes the spectral radius, strictly if $p<q$ and $s>t$.
\end{claim}

\begin{proof}[Proof of the claim]
\renewcommand{\qedsymbol}{\openbox}
Let $V^{(1)},V^{(2)}$ be the core parts of sizes $p,q$, respectively, and let $\boldsymbol b$ be the indicator vector of the neighbors of the exceptional vertex $u$ in the core.
With the retained vertex $u$ listed first, the adjacency matrix is
\[
 A(\Jn st)=
 \begin{pmatrix}
  0&\boldsymbol b^T\\
  \boldsymbol b&A_{22}
 \end{pmatrix},
 \qquad
 A_{22}=A(K_{p,q}).
\]
For $\rho>\sqrt{pq}$, write
\[
 R_{\{u\}}(\rho)=(\rho I-A_{22})^{-1},
 \qquad
 K_{\{u\}}(\rho)=
 \boldsymbol b^TR_{\{u\}}(\rho)\boldsymbol b,
\]
where the $1\times1$ effective matrix is identified with its single entry.
By Lemma~\ref{lem:core-elimination}, eliminating the core coordinates gives $\rho=K_{\{u\}}(\rho)$.
Direct inversion yields
\[
 R_{\{u\}}(\rho)
 =
 \begin{pmatrix}
 \displaystyle
 \frac{I_p}{\rho}
 +\frac{qJ_{p,p}}{\rho(\rho^2-pq)}
 &
 \displaystyle
 \frac{J_{p,q}}{\rho^2-pq}
 \\[6pt]
 \displaystyle
 \frac{J_{q,p}}{\rho^2-pq}
 &
 \displaystyle
 \frac{I_q}{\rho}
 +\frac{pJ_{q,q}}{\rho(\rho^2-pq)}
 \end{pmatrix},
\]
where each $J_{a,b}$ is an all-ones matrix.
Substitution shows that $K_{\{u\}}(\rho)$ is the right-hand side of \eqref{j:spectral-resolvent}.

For $\rho>\sqrt{pq}$, each term in $K_{\{u\}}(\rho)$ is strictly decreasing.
Indeed, the last summand can be written as $2st/(\rho^2-pq)$.
Moreover,
\[
 \lim_{\rho\to\sqrt{pq}^{+}}K_{\{u\}}(\rho)=\infty,
 \qquad
 \lim_{\rho\to\infty}K_{\{u\}}(\rho)=0.
\]
Thus the continuous decreasing function $K_{\{u\}}(\rho)$ equals $\rho$ at exactly one point.
The intersection is $\rh(\Jn st)$.

Now consider the alternative assignment in which $s$ and $t$ are swapped.
Let $K_{\{u\}}^{\mathrm{swap}}(\rho)$ be the corresponding right-hand side.
Then
\[
 K_{\{u\}}(\rho)-K_{\{u\}}^{\mathrm{swap}}(\rho)
 =
 \frac{(q-p)(s^2-t^2)}
      {\rho(\rho^2-pq)}.
\]
This difference is nonnegative, and it is strictly positive when $p<q$ and $s>t$.
The fixed-point comparison therefore gives the asserted optimal assignment and strictness.
The symmetry of $\Pprof_J(H)$ ensures that swapping the core-neighbor counts preserves admissibility.
\end{proof}


\begin{claim}\label{j:expansion}
Let $(s,t)$ range over the admissible pairs.
Then, uniformly over these pairs,
\begin{equation*}
 \rh(\Jn st)
 =
 \sqrt{pq}
 +\frac{2(s+t)^2}{(n-1)^2}
 +\Oh(n^{-3}).
\end{equation*}
\end{claim}

\begin{proof}[Proof of the claim]
\renewcommand{\qedsymbol}{\openbox}
Write $\rho=\rh(\Jn st)$ and $\sigma=s+t$.
Since $\Ttwo{n-1}$ is a proper induced subgraph of $\Jn st$, we have $\rho>\sqrt{pq}$.
The part sizes satisfy
\[
 p,q=\frac{n-1}{2}+O(1),
 \qquad
 pq=\frac{(n-1)^2}{4}+O(1).
\]
The parameters $s,t$ are bounded by constants depending only on $H$.

We first obtain crude bounds for $\rho$.
For sufficiently large $n$,
\[
 \rho^2-\sigma=\Theta(n^2),
 \qquad
 \rho^2-\sigma\ge n^2/16.
\]
Since $\rho\le n-1$, the right-hand side of \eqref{j:spectral-polynomial} satisfies
\[
 qs^2+pt^2+2\rho st=\Oh(n).
\]
Dividing by $\rho^2-\sigma=\Theta(n^2)$ gives
\[
 \rho^2-pq=\Oh(n^{-1}),
\]
and hence
\[
 \rho-\sqrt{pq}=\Oh(n^{-2}).
\]

Next, using
\[
 \sqrt{pq}-\frac{n-1}{2}
 =
 \frac{pq-(n-1)^2/4}
      {\sqrt{pq}+(n-1)/2}
 =
 O(n^{-1}),
\]
we obtain
\[
 \rho-\frac{n-1}{2}=\Oh(n^{-1}).
\]
Therefore
\begin{align*}
 qs^2+pt^2+2\rho st
 =
 \frac{n-1}{2}(s^2+t^2+2st)+\Oh(1)=
 \frac{n-1}{2}\sigma^2+\Oh(1).
\end{align*}
Also,
\[
 \rho^2-\sigma
 =
 \frac{(n-1)^2}{4}+\Oh(1).
\]
Substituting these estimates into \eqref{j:spectral-polynomial} yields
\[
 \rho^2-pq
 =
 \frac{\frac{n-1}{2}\sigma^2+\Oh(1)}
      {\frac{(n-1)^2}{4}+\Oh(1)}
 =
 \frac{2\sigma^2}{n-1}+\Oh(n^{-2}).
\]
Finally,
\[
 \rho+\sqrt{pq}=n-1+\Oh(n^{-1}),
\]
so
\[
 \rho-\sqrt{pq}
 =
 \frac{\rho^2-pq}{\rho+\sqrt{pq}}
 =
 \frac{2\sigma^2}{(n-1)^2}+\Oh(n^{-3}).
\]
\end{proof}

We now determine the optimal parameter pair. For each candidate pair, let $\sigma=s+t$.
Theorem~\ref{j:complete-core} provides the complete balanced bipartite core and bounded core-neighbor counts.
The orientation comparison in Claim~\ref{j:spectral-equation} then places the larger core neighborhood in the smaller core part.
Together with Lemma~\ref{j:frontier}, these results reduce the problem to the boundary pairs $(\FJ(t),t)$ with $t\in I_J$.
There are only finitely many such candidates.
By Claim~\ref{j:expansion}, a larger sum of the core-neighbor counts $s+t$ gives a strictly larger spectral radius for sufficiently large $n$.
Therefore any optimal pair must satisfy $t\in T_J$.

Suppose first that $n$ is odd, so $p=q=(n-1)/2$.
For a fixed sum $\sigma=s+t$, the right-hand side of \eqref{j:spectral-polynomial} becomes
\[
 p\sigma^2+2(\rho-p)st.
\]
At every relevant root, $\rho>\rho(K_{p,p})=p$.
The expression is therefore strictly increasing in $st$.
Under the ordering $s=\sigma-t\ge t$, the product $t(\sigma-t)$ strictly increases with $t$.
Thus the unique optimum among equal-sum candidates is obtained at $t=\max T_J$.

Now suppose that $n$ is even.
Let $a=(n-1)/2$, so $p=a-1/2$ and $q=a+1/2$.
Since $\sqrt{pq}=a-1/(4(n-1))+O(n^{-3})$, Claim~\ref{j:expansion} gives
\[
 \rho-a
 =
 -\frac{1}{4(n-1)}
 +\frac{2\sigma^2}{(n-1)^2}
 +\Oh(n^{-3})
 <0
\]
at every relevant root for sufficiently large $n$.

Let $\tau=s-t\ge0$.
Since $s+t=\sigma$, we have
\[
 s^2+t^2=\frac{\sigma^2+\tau^2}{2},
 \qquad
 st=\frac{\sigma^2-\tau^2}{4}.
\]
Thus the right-hand side of
\eqref{j:spectral-polynomial} can be written as
\[
 a\sigma^2+\frac{\sigma\tau}{2}
 +\frac{\rho-a}{2}(\sigma^2-\tau^2).
\]
For fixed $\sigma$ and $\rho$, its derivative with respect
to $\tau$ is $\frac{\sigma}{2}-(\rho-a)\tau>0.$
Hence it strictly increases with $s-t$.
For a fixed sum $s+t$, this difference is maximized by
taking $t=\min T_J$.

By Claim~\ref{j:spectral-equation}, these strict comparisons at the relevant roots give strict comparisons of the corresponding spectral radii.
This proves the parameter rule and uniqueness.
\end{proof}

The following corollary identifies the symmetric, book-type extremal graph directly from the profile of $H$.

\begin{corollary}\label{j:diagonal}
Let $m=\dst(H)-2\ge1$.
For sufficiently large $n$, the unique extremal graph is $\Jn mm$ if and only if
$\max\{\alpha,\beta\}\ge m+1,$ for every $(\alpha,\beta)\in\Pprof_J(H).$
\end{corollary}

\begin{proof}
By Lemma~\ref{j:frontier}, the displayed condition is exactly the condition that $\Jn mm$ is $H$-free.
Every admissible pair satisfies $s,t\le m$.
Thus, if $(m,m)$ is admissible, it uniquely maximizes $s+t$, and hence the spectral radius.
Conversely, an extremal graph must be $H$-free, which gives the displayed condition.
\end{proof}

\subsection{The case $\dst(H)=2$ with a triangle:
the subdivided Tur\'an graph}\label{sec:trianglecase}

In this subsection, we prove that the extremal graph itself is triangle-free.
By Lemma~\ref{prelim:triangle}, this structural property forces the graph to be $\ST{n-1}$.

\begin{theorem}\label{s:structural}
If $\dst(H)=2$ and $H$ contains a triangle, then
\[
 \SPEX_3(n,H)=\{\ST{n-1}\}
\]
for sufficiently large $n$.
\end{theorem}

\begin{proof}
Choose $u_1v_1\in\Estar(H)$ with $d_H(u_1)=2$, and write $N_H(u_1)=\{v_1,v_2\}$.
Lemma~\ref{prelim:critical} gives $v_1v_2\in E(H)$. Recall $n_0=|V(H)|$.

Suppose that $G$ contains a triangle.
Since $G-u$ is bipartite, this triangle is $uu_1^{(1)}v_1^{(2)}$, where $u_1^{(1)}\in V^{(1)}$, $v_1^{(2)}\in V^{(2)}$, and $u_1^{(1)}v_1^{(2)}\in E(B)$.
By Lemmas~\ref{prelim:rough} and~\ref{prelim:rectangle},
\[
 \min\{d_B(u_1^{(1)}),d_B(v_1^{(2)})\}
 \le
 n_0\sqrt{\bar e(B)+1}+1,
 \qquad
 \bar e(B)<3n.
\]
Consequently, one of $u_1^{(1)},v_1^{(2)}$, say $v$, satisfies $d_G(v)=\Oh(\sqrt n)$.
Hence
\[
 x_v
 \le
 \frac{d_G(v)\pxmax}{\rho}
 =
 \Oh(n^{-1}).
\]
Minimality of $x_u$ gives the same bound for $x_u$.

Deleting $u,v$ and applying Lemma~\ref{prelim:delete} therefore gives
\[
 \rh(G-\{u,v\})^2
 >
 e(\Ttwo{n-1})-2-\Oh(1).
\]
But $G-\{u,v\}$ is bipartite of order $n-2$, so
\[
 \rh(G-\{u,v\})^2\le e(\Ttwo{n-2}).
\]
This is impossible because
\[
 e(\Ttwo{n-1})-e(\Ttwo{n-2})
 =
 \left\lfloor\frac{n-1}{2}\right\rfloor
 \longrightarrow\infty.
\]

Thus $G$ is triangle-free.
The graph $\ST{n-1}$ is admissible because $H$ contains a triangle.
Lemma~\ref{prelim:triangle} proves the asserted maximum and uniqueness.
\end{proof}

\subsection{The case $\dst(H)=2$ without a triangle:
the $L$-family}
\label{sec:Lparameters}

In this subsection,  we show that the extremal graph belongs to the \(L\)-family and determine it exactly.


The following proposition establishes the complete core and the bounded exceptional neighborhoods.


\begin{proposition}\label{l:structural}
Let $H$ be triangle-free with $\dst(H)=2$, $G\in\SPEX_3(n,H)$, and $u$ be a vertex with $x_u=\pxmin$.
For sufficiently large $n$, there exists a vertex $v$ such that the following hold:
\begin{enumerate}[label=(\roman*),leftmargin=*]
\item
The vertices $u,v$ are adjacent, and $B'=G-\{u,v\}$ is bipartite.

\item
The sets $U=N_G(u)\setminus\{v\},~V=N_G(v)\setminus\{u\}$ lie in the same part of $B'$ and satisfy $1\le|U|,|V|<n_0.$
In particular, both exceptional neighborhoods are nonempty and have sizes $\Oh(1)$.

\item
The remaining graph is a complete balanced bipartite core: $B'\cong\Ttwo{n-2}.$
\end{enumerate}
\end{proposition}

\begin{proof}
Choose a critical edge $u_1v_1$ with $d_H(u_1)=2$, and let $N_H(u_1)=\{v_1,v_2\}$.
Since $H$ is triangle-free, $v_1v_2\notin E(H)$.
Recall that $B=G-u$ is bipartite with parts $V^{(1)},V^{(2)}$.

$(i).$ We first identify the second exceptional vertex $v$ and obtain preliminary estimates for the remaining core. Choose
$u_1^{(1)}\in N_G(u)\cap V^{(1)},~v_1^{(2)}\in N_G(u)\cap V^{(2)}.$
Lemma~\ref{prelim:rectangle}, which requires no edge $u_1^{(1)}v_1^{(2)}$ in the present case, together with Lemma~\ref{prelim:rough}, provides a vertex $v\in\{u_1^{(1)},v_1^{(2)}\}\subseteq N_G(u)$ satisfying $d_G(v)=\Oh(\sqrt n)$.
Thus $uv\in E(G)$.
The eigenvalue equation, together with $\pxmax\le3/\sqrt n$ and $\rho\ge n/3$, gives
\[
 x_v
 \le
 \frac{d_G(v)\pxmax}{\rho}
 =
 \Oh(n^{-1}).
\]
By minimality of $x_u$,
\[
 x_u\le x_v=\Oh(n^{-1}).
\]

The graph $B'=G-\{u,v\}$ is bipartite because it is a subgraph of $B$.
By Lemmas~\ref{prelim:lower}, \ref{prelim:delete}, and~\ref{prelim:bipartite},
\[
 e(B')\ge\rh(B')^2
 >
 e(\Ttwo{n-2})-\Oh(1).
\]
Thus $e(\Ttwo{n-2})-e(B')=\Oh(1)$.
If $p,q$ are the part sizes of $B'$, then
\[
 \bigl(e(\Ttwo{n-2})-pq\bigr)+\bar e(B')=\Oh(1).
\]
Both summands are nonnegative.
Consequently,
\[
 \bar e(B')=\Oh(1),
 \qquad
 p,q=\frac n2+\Oh(1).
\]
Every vertex of $B'$ therefore has degree $n/2-\Oh(1)$ in $B'$.

$(ii).$ We show that the exceptional neighborhoods are nonempty, lie in the same core part, and have bounded sizes. After interchanging the original parts and relabeling the chosen core vertices if necessary, assume that $v=v_1^{(2)}$.
The parts of $B'$ are then $V^{(1)}$ and $V^{(2)}\setminus\{v\}$, and $V=N_G(v)\setminus\{u\}\subseteq V^{(1)}.$
If $u$ had core neighbors in both parts of $B'$, Lemma~\ref{prelim:rectangle}, applied to $G-v$, would give a core neighbor of $u$ with degree $\Oh(1)$ in $B'$, because $\bar e(B')=\Oh(1)$.
This contradicts the degree estimate obtained in (i).

The core neighborhood of $u$ cannot be empty or lie only in $V^{(2)}\setminus\{v\}$, since either possibility would make $G$ bipartite.
Hence $U=N_G(u)\setminus\{v\}\subseteq V^{(1)},$ and $U\ne\varnothing.$ Similarly, $V\ne\varnothing$.

We claim that $|V|<n_0$.
Suppose instead that $|V|\ge n_0$.
Since $u_1^{(1)}\in U$, map
\[
 u_1\mapsto u,
 \qquad
 v_1\mapsto v,
 \qquad
 v_2\mapsto u_1^{(1)}.
\]
In a bipartition of $H-u_1v_1$, the vertices $u_1,v_1$ have the same color, while $v_2$ has the opposite color.
Map every other vertex of the color of $v_2$ into $V\setminus\{u_1^{(1)}\}$.
There are enough available vertices.
All vertices of this color other than $v_2$ are mapped into $V$, and $v_1v_2\notin E(H)$, so every required edge incident with $v_1$ is present.

Map the remaining vertices into $V^{(2)}\setminus\{v\}$.
The previously chosen core vertices exclude only $\Oh(1)$ vertices through missing adjacencies, because $\bar e(B')=\Oh(1)$.
Since $|V^{(2)}\setminus\{v\}|=\frac n2+\Oh(1),$
these remaining choices can be made injectively for sufficiently large $n$.
This embeds $H$ in $G$, a contradiction.
Thus $|V|<n_0$.
Interchanging $u,v$ gives $|U|<n_0$.

$(iii).$ We sharpen the deletion estimate to prove that the core is exactly complete and balanced. By $(ii)$, we have
$d_G(u)=|U|+1\le n_0,
 d_G(v)=|V|+1\le n_0.$
The eigenvalue equations therefore imply
$x_u,x_v
 \le
 \frac{n_0\pxmax}{\rho}
 =
 \Oh(n^{-3/2}).$
Applying Lemma~\ref{prelim:delete} again yields
\[
 \rh(B')^2
 >
 e(\Ttwo{n-2})-\Oh(n^{-1}).
\]
Together with Lemma~\ref{prelim:bipartite}, this gives
\[
 e(\Ttwo{n-2})-\Oh(n^{-1})
 <
 \rh(B')^2
 \le e(B')
 \le pq
 \le e(\Ttwo{n-2}).
\]
Since $e(B')$ is an integer, sufficiently large $n$ forces
\[
 e(B')=pq=e(\Ttwo{n-2}).
\]
Hence $B'\cong\Ttwo{n-2}$, as required.
\end{proof}

The next containment criterion allows us to modify the
exceptional neighborhoods while preserving $H$-freeness.
To state it, let $G'$ consist of a complete bipartite core
with parts $V^{(1)},V^{(2)}$ (satisfying $|V^{(1)}|,|V^{(2)}|>n_0$)
and two adjacent vertices $u,v$ with nonempty core
neighborhoods $U,V\subseteq V^{(1)}$, respectively.

\begin{lemma}\label{l:arbitrary-containment}
If $H$ is a triangle-free three-chromatic color-critical graph, then $G'$ contains $H$ if and only if there is an oriented critical edge $u_1v_1\in\Estar(H)$ satisfying
\[
 d_H(u_1)-1\le|U|,
 \qquad
 d_H(v_1)-1\le|V|,
 \qquad
 d_H(u_1)+d_H(v_1)-2\le|U\cup V|.
\]
\end{lemma}

\begin{proof}
Since $G'-uv$ is bipartite, every copy of $H$ must use $uv$, whose preimage is a critical edge $u_1v_1$.
Orient $u_1v_1$ so that $u_1,v_1$ map to $u,v$.
Triangle-freeness makes the residual neighborhoods
\[
 S_{u_1}=N_H(u_1)\setminus\{v_1\},
 \qquad
 S_{v_1}=N_H(v_1)\setminus\{u_1\}
\]
disjoint.
Their images lie in $U,V$, respectively, so injectivity gives the three necessary inequalities.

Conversely, suppose an oriented critical edge satisfies these inequalities.
Construct an auxiliary bipartite graph with left classes indexed by $S_{u_1},S_{v_1}$ and right vertex set $U\cup V$.
Give each vertex of $S_{u_1}$ neighborhood $U$, and each vertex of $S_{v_1}$ neighborhood $V$.
For a nonempty subset $S$ of the left vertex set, its neighborhood is $U$, $V$, or $U\cup V$, according as $S$ meets only the first class, only the second class, or both.
The three inequalities imply $|N(S)|\ge|S|$ in every case.
Hall's marriage theorem supplies a matching saturating the left vertex set.

Map $u_1,v_1$ to $u,v$ and their residual neighborhoods according to the matching.
In a bipartition of $H-u_1v_1$, the vertices $u_1,v_1$ have the same color, so all remaining vertices embed into unused vertices of the corresponding core parts.
Both parts have at least $n_0$ vertices.
This gives a copy of $H$ in $G'$.
\end{proof}

The following theorem identifies the $L$-family.
Its proof establishes the nesting and orientation properties through two claims.

\begin{theorem}\label{l:family}
Let $H$ be a fixed connected triangle-free three-chromatic color-critical graph with $\dst(H)=2$.
If $n$ is sufficiently large and $G\in\SPEX_3(n,H)$, then $G\cong\Ln st$, for some $1\le t\le s<n_0.$
In particular, the two adjacent vertices $u,v$ have nonempty nested neighborhoods in a smallest part of $\Ttwo{n-2}$.
\end{theorem}

\begin{proof}
By Proposition~\ref{l:structural}, deleting two adjacent vertices $u,v$ leaves a core
$B'\cong\Ttwo{n-2}$, and the nonempty sets $U=N_G(u)\setminus\{v\}$ and $V=N_G(v)\setminus\{u\}$
lie in a single core part and have sizes less than $n_0$.
Label the core parts so that $U,V\subseteq V^{(1)}$.
We first prove that these neighborhoods are nested.

\begin{claim}\label{l:nested}
The sets $U,V$ are nested.
\end{claim}

\begin{proof}[Proof of the claim]
\renewcommand{\qedsymbol}{\openbox}
Suppose otherwise.
Relabel $u,v$ so that $x_u\ge x_v$, and choose
$u_1^{(1)}\in U\setminus V,\ u_2^{(1)}\in V\setminus U.$
Because the core is complete bipartite, the eigenvalue equations at $u_1^{(1)},u_2^{(1)}$ give
\[
 \rho(x_{u_1^{(1)}}-x_{u_2^{(1)}})
 =
 x_u-x_v
 \ge0.
\]
Replace the edge $vu_2^{(1)}$ by $vu_1^{(1)}$.

The individual neighborhood sizes are unchanged, while the size of their union decreases by one.
Lemma~\ref{l:arbitrary-containment} therefore shows that the new graph remains $H$-free.
It is non-bipartite because it contains the triangle $uvu_1^{(1)}$.

The Rayleigh gain is
$2x_v(x_{u_1^{(1)}}-x_{u_2^{(1)}})\ge0.$
If the gain is zero, the old Perron vector is not an eigenvector of the new graph.
The Rayleigh principle then still gives a strict increase in spectral radius.
Either possibility contradicts extremality.
\end{proof}

Relabel $u,v$ so that $V\subseteq U$, then $|U|=s\ge t=|V|\ge1$.
Let $p=\left\lfloor(n-2)/2\right\rfloor$ and $q=\left\lceil(n-2)/2\right\rceil$.
Since $s<n_0$ is bounded, $s\le p$ for sufficiently large $n$.
We next determine which part of the core contains these neighborhood sets $V$ and $U$.

\begin{claim}\label{l:orientation}
For fixed neighborhood sizes $s,t$, placing both nested neighborhoods in a smallest part of $B'\cong K_{p,q}$ maximizes the spectral radius.
The preference is strict when $p<q$.
\end{claim}

\begin{proof}[Proof of the claim]
\renewcommand{\qedsymbol}{\openbox}
Write $V^{(1)},V^{(2)}$ for the core parts, with the core
neighborhoods lying in $V^{(1)}$.
Thus $\{|V^{(1)}|,|V^{(2)}|\}=\{p,q\}$.
Let $W$ be the $|V^{(1)}|\times2$ incidence matrix whose
columns are the indicator vectors of $U$ and $V$.
Let $J_{|V^{(1)}|,|V^{(2)}|}$ denote the
$|V^{(1)}|\times|V^{(2)}|$ all-ones matrix, and let
\[
 A_{11}=
 \begin{pmatrix}
  0&1\\
  1&0
 \end{pmatrix},
 ~
 A_{12}=\begin{pmatrix}W^T&0\end{pmatrix},
 ~
 A_{22}=A(K_{|V^{(1)}|,|V^{(2)}|})=
 \begin{pmatrix}
  0&J_{|V^{(1)}|,|V^{(2)}|}\\
  J_{|V^{(2)}|,|V^{(1)}|}&0
 \end{pmatrix}.
\]
With the vertices ordered as $u,v,V^{(1)},V^{(2)}$, we have
\[
 A(G)=
 \begin{pmatrix}
  A_{11}&W^T&0\\
  W&0&J_{|V^{(1)}|,|V^{(2)}|}\\
  0&J_{|V^{(2)}|,|V^{(1)}|}&0
 \end{pmatrix}
 =
 \begin{pmatrix}
  A_{11}&A_{12}\\
  A_{12}^T&A_{22}
 \end{pmatrix}.
\]
Since $K_{p,q}\cong\Ttwo{n-2}$ is a proper induced subgraph
of the connected graph $G$, we have $\rho>\sqrt{pq}$.

Write
$\boldsymbol z=\px|_{\{u,v\}}$ and
$\boldsymbol y=\px|_{V^{(1)}\cup V^{(2)}}$.
The Perron equations become
\[
 A_{11}\boldsymbol z+A_{12}\boldsymbol y
 =\rho\boldsymbol z,\qquad
 A_{12}^T\boldsymbol z+A_{22}\boldsymbol y
 =\rho\boldsymbol y.
\]
Let $M_\rho
 :=
 \rho I_{|V^{(1)}|+|V^{(2)}|}-A_{22}.$
Because $\rho>\rh(A_{22})=\sqrt{pq}$, the matrix $M_\rho$
is invertible.
The second Perron equation therefore gives
$\boldsymbol y
 =
 R_{\{u,v\}}(\rho)A_{12}^T\boldsymbol z$ and
$R_{\{u,v\}}(\rho)=M_\rho^{-1}.$
Substitution into the first equation, as in
Lemma~\ref{lem:core-elimination}, yields $K_{\{u,v\}}(\rho)\boldsymbol z=\rho\boldsymbol z,$
where
\[
 \begin{aligned}
 K_{\{u,v\}}(\rho)
 &=
 A_{11}+A_{12}R_{\{u,v\}}(\rho)A_{12}^T=
 A_{11}
 +
 \begin{pmatrix}W^T&0\end{pmatrix}
 R_{\{u,v\}}(\rho)
 \begin{pmatrix}W\\0\end{pmatrix}\\
 &=
 A_{11}
 +
 W^T
 \bigl[R_{\{u,v\}}(\rho)\bigr]_{V^{(1)}V^{(1)}}
 W.
 \end{aligned}
\]
$\bigl[R_{\{u,v\}}(\rho)\bigr]_{V^{(1)}V^{(1)}}$ is a block of $R_{\{u,v\}}(\rho)$ indexed by
$V^{(1)}$. Using the Schur-complement notation of
Definition~\ref{def:schur-complement}, after simultaneously interchanging
the two block rows and the corresponding block columns, we obtain
\[
 (\rho I_n-A(G))/M_\rho
 =
 \rho I_2-A_{11}-A_{12}M_\rho^{-1}A_{12}^T
 =
 \rho I_2-K_{\{u,v\}}(\rho).
\]
In particular, Lemma~\ref{lem:schur-determinant} gives
$\det(\rho I_n-A(G))
 =
 \det(M_\rho)
 \det(\rho I_2-K_{\{u,v\}}(\rho)).$
Now compute the required block
$\bigl[R_{\{u,v\}}(\rho)\bigr]_{V^{(1)}V^{(1)}}
 =
 \bigl[M_\rho^{-1}\bigr]_{V^{(1)}V^{(1)}}$.
By the definition of $A_{22}$,
\[
 M_\rho
 =
 \begin{pmatrix}
  \rho I_{|V^{(1)}|}
  &-J_{|V^{(1)}|,|V^{(2)}|}\\
  -J_{|V^{(2)}|,|V^{(1)}|}
  &\rho I_{|V^{(2)}|}
 \end{pmatrix}.
\]
By the block inverse formula and the notation of
Definition~\ref{def:schur-complement}, the block
$\bigl[M_\rho^{-1}\bigr]_{V^{(1)}V^{(1)}}$
is the inverse of the Schur complement of the lower-right block
$\rho I_{|V^{(2)}|}$ in $M_\rho$. Thus
\[
 \begin{aligned}
 \bigl[R_{\{u,v\}}(\rho)\bigr]_{V^{(1)}V^{(1)}}
 &=
 \bigl(M_\rho/(\rho I_{|V^{(2)}|})\bigr)^{-1}=
 \left(
 \rho I_{|V^{(1)}|}
 -
 J_{|V^{(1)}|,|V^{(2)}|}
 (\rho I_{|V^{(2)}|})^{-1}
 J_{|V^{(2)}|,|V^{(1)}|}
 \right)^{-1}\\
 &=
 \left(
 \rho I_{|V^{(1)}|}
 -
 \frac{|V^{(2)}|}{\rho}
 J_{|V^{(1)}|}J_{|V^{(1)}|}^T
 \right)^{-1}.
 \end{aligned}
\]
Here $J_{|V^{(1)}|}$ is the all-ones column vector of
length $|V^{(1)}|$. Applying the rank-one inverse formula and using
$|V^{(1)}||V^{(2)}|=pq$, we obtain
\[
 \begin{aligned}
 \bigl[R_{\{u,v\}}(\rho)\bigr]_{V^{(1)}V^{(1)}}
 &=
 \frac{I_{|V^{(1)}|}}{\rho}
 +
 \frac{|V^{(2)}|/\rho^3}
      {1-|V^{(1)}||V^{(2)}|/\rho^2}
 J_{|V^{(1)}|}J_{|V^{(1)}|}^T\\
 &=
 \frac{I_{|V^{(1)}|}}{\rho}
 +
 \frac{|V^{(2)}|}{\rho(\rho^2-pq)}
 J_{|V^{(1)}|}J_{|V^{(1)}|}^T.
 \end{aligned}
\]
The denominator is positive because $\rho>\sqrt{pq}$.

The nesting $V\subseteq U$ gives
$
 W^TW=
 \begin{pmatrix}
  s&t\\
  t&t
 \end{pmatrix},~
 W^TJ_{|V^{(1)}|}
 =
 \boldsymbol w
 =
 \begin{pmatrix}s\\t\end{pmatrix}.
$
Consequently,
$
 W^TJ_{|V^{(1)}|}J_{|V^{(1)}|}^TW
 =
 \boldsymbol w\boldsymbol w^T
 =
 \begin{pmatrix}
  s^2&st\\
  st&t^2
 \end{pmatrix}.
$
To compare the two possible placements, write
\[
 K_{\{u,v\}}(\rho;c)
 :=
 A_{11}
 +
 \frac1\rho
 \begin{pmatrix}
  s&t\\
  t&t
 \end{pmatrix}
 +
 \frac{c}{\rho(\rho^2-pq)}
 \boldsymbol w\boldsymbol w^T,
 \qquad c\in\{p,q\}.
\]
Substituting the computed resolvent block into the
effective matrix gives
\begin{equation}\label{l:effective-matrix}
 \begin{aligned}
 K_{\{u,v\}}(\rho)
 &=K_{\{u,v\}}(\rho;|V^{(2)}|)=
 A_{11}
 +
 \frac1\rho
 \begin{pmatrix}
  s&t\\
  t&t
 \end{pmatrix}
 +
 \frac{|V^{(2)}|}{\rho(\rho^2-pq)}
 \boldsymbol w\boldsymbol w^T,\\
 \rho&=\rh(K_{\{u,v\}}(\rho)).
 \end{aligned}
\end{equation}
The spectral-radius identity follows from
Lemma~\ref{lem:core-elimination}: the vector
$\boldsymbol z$ is positive and $K_{\{u,v\}}(\rho)$
is a positive matrix.

For fixed $c\in\{p,q\}$, every entry of
$K_{\{u,v\}}(\rho;c)$ strictly decreases as $\rho$
increases on $(\sqrt{pq},\infty)$.
Thus $\rh(K_{\{u,v\}}(\rho;c))$ is continuous and
strictly decreasing on this interval.
Moreover,
\[
 \lim_{\rho\to\sqrt{pq}^{+}}
 \rh(K_{\{u,v\}}(\rho;c))=\infty,
 \qquad
 \lim_{\rho\to\infty}
 \rh(K_{\{u,v\}}(\rho;c))=1.
\]
Hence $\rho=\rh(K_{\{u,v\}}(\rho;c))$ has a unique
solution on $(\sqrt{pq},\infty)$, namely the spectral
radius of the graph with $|V^{(2)}|=c$.
For fixed $pq$, increasing $c$ strictly increases every
entry of $K_{\{u,v\}}(\rho;c)$, so strict Perron comparison
shows that this solution also strictly increases.

Finally, placing the neighborhoods in the part of size
$p$ corresponds to $|V^{(2)}|=q$, whereas placing the
neighborhoods in the part of size $q$ corresponds to
$|V^{(2)}|=p$.
Since $p\le q$, for every $\rho>\sqrt{pq}$ we have
\[
 K_{\{u,v\}}(\rho;q)-K_{\{u,v\}}(\rho;p)
 =
 \frac{q-p}{\rho(\rho^2-pq)}
 \boldsymbol w\boldsymbol w^T
 \ge0
\]
entrywise.
Thus $\rh(K_{\{u,v\}}(\rho;q))\ge\rh(K_{\{u,v\}}(\rho;p)),$
with strict inequality when $p<q$.
The fixed-point comparison proves the claimed preference.
For sufficiently large $n$,
Lemma~\ref{l:arbitrary-containment} shows that changing
the part containing the neighborhoods preserves
$H$-freeness, since $|U|$, $|V|$, and $|U\cup V|$
remain unchanged.
\end{proof}

By Claim~\ref{l:orientation}, spectral extremality places the two neighborhoods in a smallest part; when the parts have equal sizes, the two choices give isomorphic graphs.
Together with Claim~\ref{l:nested}, these are precisely the defining properties of $\Ln st$.
\end{proof}

We now identify the unique nested parameter pair.
The endpoint degrees of the critical edges give a containment criterion, and a spectral expansion will yield a parameter rule independent of the parity of $n$.

For each critical edge $u_1v_1$, let
\begin{equation*}
 \alpha_{u_1v_1}
 =
 \min\{d_H(u_1)-1,d_H(v_1)-1\},
 \qquad
 \beta_{u_1v_1}
 =
 d_H(u_1)+d_H(v_1)-2.
\end{equation*}
Define $\Pprof_L(H)
=
\{(\alpha_{u_1v_1},\beta_{u_1v_1}):
u_1v_1\in\Estar(H)\}$
to be \emph{the critical-edge degree profile set}.
This definition applies to every three-chromatic color-critical graph, including disconnected graphs.
Let $\FL(t)
=
\min\{\beta-1:
(\alpha,\beta)\in\Pprof_L(H),\ \alpha\le t\}$
and $k=\min\{\beta-1:(\alpha,\beta)\in\Pprof_L(H)\}$.
Since $\dst(H)=2$, some critical edge satisfies
$\alpha_{u_1v_1}=1$.
Hence $\FL(t)$ is well defined for every positive integer $t$.

The following lemma gives the two correction terms needed to compare the possible parameter pairs.
Its finite range ensures that all remainder estimates are uniform.

\begin{lemma}\label{l:expansion}
Let
$
 p=\left\lfloor\frac{n-2}{2}\right\rfloor,
 q=\left\lceil\frac{n-2}{2}\right\rceil.
$
For sufficiently large $n$, uniformly over the integer pairs
$1\le t\le s\le n_0-2$, we have
\begin{equation}\label{l:asymptotic-radius}
 \rh(\Ln st)
 =
 \sqrt{pq}
 +\frac{q(s^2+t^2)}{2(pq)^{3/2}}
 +\frac{qst}{(pq)^2}
 +\Oh(n^{-4}).
\end{equation}
\end{lemma}

\begin{proof}
Write $\rho=\rh(\Ln st)$.
Let
$
 A_{11}=
 \begin{pmatrix}
  0&1\\
  1&0
 \end{pmatrix},
 B_0=
 \begin{pmatrix}
  s&t\\
  t&t
 \end{pmatrix},
 \boldsymbol w=
 \begin{pmatrix}s\\t\end{pmatrix}.
$
The block calculation in Claim~\ref{l:orientation} in the proof of Theorem~\ref{l:family} applies to every graph $\Ln st$.
Thus, by \eqref{l:effective-matrix} and Lemma~\ref{lem:core-elimination},
\[
 \rho\boldsymbol z
 =
 K_{\{u,v\}}(\rho)\boldsymbol z
 =
 K_{\{u,v\}}(\rho;q)\boldsymbol z
 =
 \left(
 A_{11}+\frac{B_0}{\rho}
 +\frac{q}{\rho(\rho^2-pq)}
 \boldsymbol w\boldsymbol w^T
 \right)\boldsymbol z,
\]
where $\boldsymbol z$ is the restriction of the positive Perron vector to the two exceptional vertices.

Set
$
 C_\rho
 =
 \rho I-A_{11}-\frac{B_0}{\rho}.
$
Since $\rho>\sqrt{pq}$ and $s,t$ are bounded, $C_\rho$ is positive definite for sufficiently large $n$.
Rearranging the eigenvalue equation gives
$
 C_\rho\boldsymbol z
 =
 \frac{q}{\rho(\rho^2-pq)}
 \boldsymbol w(\boldsymbol w^T\boldsymbol z).
$
Multiplying by $C_\rho^{-1}$ and then by $\boldsymbol w^T$, and canceling the positive scalar $\boldsymbol w^T\boldsymbol z$, yields
\begin{equation}\label{l:rank-one-equation}
 1
 =
 \frac{q}{\rho(\rho^2-pq)}
 \boldsymbol w^TC_\rho^{-1}\boldsymbol w.
\end{equation}

We now invert $C_\rho$ to the required order.
Let
$
 E_\rho
 =
 \frac{A_{11}}{\rho}
 +
 \frac{B_0}{\rho^2},~
 C_\rho=\rho(I-E_\rho).
$
Because $\rho=\Theta(n)$,
$
 \|E_\rho\|=\Oh(n^{-1})\le\frac12
$
for large $n$.
The Neumann series identity in Lemma~\ref{lem:resolvent-neumann} gives
\[
 C_\rho^{-1}
 =
 \frac1\rho
 \bigl(
 I+E_\rho+
 E_\rho^2(I-E_\rho)^{-1}
 \bigr)
 =
 \frac I\rho+\frac{A_{11}}{\rho^2}+Q_\rho,
\]
where the remainder matrix is
$
 Q_\rho
 =
 \frac{B_0}{\rho^3}
 +
 \frac1\rho
 E_\rho^2(I-E_\rho)^{-1}.
$
The norm estimate in Lemma~\ref{lem:resolvent-neumann} gives $\|(I-E_\rho)^{-1}\|\le2$, so $\|Q_\rho\|=\Oh(n^{-3}).$

Using
$
 \boldsymbol w^T\boldsymbol w=s^2+t^2,~
 \boldsymbol w^TA_{11}\boldsymbol w=2st,
$
we obtain
\[
 \boldsymbol w^TC_\rho^{-1}\boldsymbol w
 =
 \frac{s^2+t^2}{\rho}
 +
 \frac{2st}{\rho^2}
 +
 \Oh(n^{-3}).
\]
Substituting into \eqref{l:rank-one-equation} and multiplying by $\rho^2(\rho^2-pq)$ gives
\begin{equation}\label{l:expanded-equation}
 \rho^2(\rho^2-pq)
 =
 q(s^2+t^2)+\frac{2qst}{\rho}
 +\Oh(n^{-1}).
\end{equation}
The resulting error is $\Oh(n^{-1})$, since the $\Oh(n^{-3})$ remainder is multiplied by $q\rho=\Theta(n^2)$.

Equation~\eqref{l:expanded-equation} first implies
$
 \rho^2-pq=\Oh(n^{-1}),~
 \varepsilon:=\rho-\sqrt{pq}=\Oh(n^{-2}).
$
The exact expansion of its left-hand side is
\begin{align*}
 \rho^2(\rho^2-pq)
 &=
 (\sqrt{pq}+\varepsilon)^2
 \bigl(2\sqrt{pq}\,\varepsilon+\varepsilon^2\bigr)=
 2(pq)^{3/2}\varepsilon
 +5pq\,\varepsilon^2
 +4\sqrt{pq}\,\varepsilon^3
 +\varepsilon^4.
\end{align*}
Since $\varepsilon=\Oh(n^{-2})$, the last three terms have orders $\Oh(n^{-2})$, $\Oh(n^{-5})$, and $\Oh(n^{-8})$, respectively.

Also,
$
 \frac1\rho-\frac1{\sqrt{pq}}
 =
 -\frac{\varepsilon}{\rho\sqrt{pq}}
 =
 \Oh(n^{-4}),
$
so
$
 \frac{2qst}{\rho}
 =
 \frac{2qst}{\sqrt{pq}}+\Oh(n^{-3}).
$
Substituting these estimates into \eqref{l:expanded-equation} yields
\[
 2(pq)^{3/2}\varepsilon
 =
 q(s^2+t^2)
 +
 \frac{2qst}{\sqrt{pq}}
 +
 \Oh(n^{-1}).
\]
Dividing by $2(pq)^{3/2}$ gives
\[
 \varepsilon
 =
 \frac{q(s^2+t^2)}{2(pq)^{3/2}}
 +
 \frac{qst}{(pq)^2}
 +
 \Oh(n^{-4}),
\]
as required.
\end{proof}

We next use the spectral expansion to select the unique optimal pair.
Choose $t_L\in\{1,\ldots,k\}$ by lexicographic maximization:
\begin{equation}\label{l:selected-parameters}
 \bigl(\FL(t_L)^2+t_L^2,\ \FL(t_L)t_L\bigr)
 =
 \max_{1\le t\le k}^{\mathrm{lex}}
 \bigl(\FL(t)^2+t^2,\ \FL(t)t\bigr),
\end{equation}
and set $s_L=\FL(t_L)$.
Thus the first coordinate is maximized first, and the second is used only to break a tie.

The following theorem proves that this rule determines the exact parameters.
The containment and boundary reductions are proved within its proof.

\begin{theorem}\label{l:parameters}
Suppose that $\dst(H)=2$ and $H$ is triangle-free.
For sufficiently large $n$,
\[
 \SPEX_3(n,H)=\{\Ln{s_L}{t_L}\}.
\]
The ordered pair $(s_L,t_L)$ is uniquely determined by \eqref{l:selected-parameters} and is independent of the parity of $n$.
\end{theorem}

\begin{proof}
By Theorem~\ref{l:family}, every spectral extremal graph has the form $\Ln st$ with bounded integers $s\ge t\ge1$.
We first determine which pairs are admissible.

\begin{claim}\label{l:nested-containment}
Let $s\ge t\ge1$ be fixed integers.
Then, for sufficiently large $n$, $H\subseteq\Ln st$ if and only if there exists a critical edge $u_1v_1\in\Estar(H)$ whose profile satisfies
\begin{equation}\label{l:containment}
 \alpha_{u_1v_1}\le t
 \quad\text{and}\quad
 \beta_{u_1v_1}\le s.
\end{equation}
\end{claim}

\begin{proof}[Proof of the claim]
\renewcommand{\qedsymbol}{\openbox}
In the nested graph, the union of the core neighborhoods has size $s$ and the smaller neighborhood has size $t$.
Necessity follows from Lemma~\ref{l:arbitrary-containment}.

Conversely, orient $u_1v_1$ so that $d_H(v_1)-1=\alpha_{u_1v_1}$, and map $v_1$ to the exceptional vertex with the smaller neighborhood.
The other residual degree is at most $\beta_{u_1v_1}\le s$.
Together with the smaller-neighborhood bound and the union bound, these are exactly the conditions in Lemma~\ref{l:arbitrary-containment}.
\end{proof}

\begin{claim}\label{l:frontier}
Let $s\ge t\ge1$.
Then $(s,t)$ is admissible if and only if $s\le\FL(t)$.
The function $\FL$ is nonincreasing, satisfies $\FL(k)=k$, and the feasible values of $t$ are exactly $1,\ldots,k$.
For each such $t$, the largest spectral radius is attained at $s=\FL(t)$.
Thus the boundary candidates are
\begin{equation}\label{l:boundary-candidates}
 (\FL(t),t),
 \qquad
 1\le t\le k.
\end{equation}
\end{claim}

\begin{proof}[Proof of the claim]
\renewcommand{\qedsymbol}{\openbox}
The admissibility condition follows directly from Claim~\ref{l:nested-containment}, equivalently from \eqref{l:containment}.
Monotonicity follows from the definition of $\FL$.

Choose a critical edge attaining $\beta_{u_1v_1}-1=k$.
Then
\[
 \alpha_{u_1v_1}
 \le
 \frac{\beta_{u_1v_1}}2
 =
 \frac{k+1}{2}
 \le k.
\]
Therefore $\FL(k)=k$.
For $t\le k$, every value in the minimum is at least $k$, so $\FL(t)\ge k\ge t$.
For $t>k$, monotonicity and the global minimum give $\FL(t)=k<t$, so no admissible pair can have $t>k$.

Increasing $s$ while keeping $t$ fixed adds edges to a connected graph and strictly increases its spectral radius.
\end{proof}

By Claim~\ref{l:frontier}, the optimization is over the finite boundary list~\eqref{l:boundary-candidates}.
We check that the expansion in Lemma~\ref{l:expansion} applies uniformly to all its candidates.
Choose a critical edge $u_1v_1$ with $d_H(u_1)=2$.
Then $\alpha_{u_1v_1}=1$ and $\beta_{u_1v_1}=d_H(v_1)$, so
\[
 \FL(1)\le\beta_{u_1v_1}-1=d_H(v_1)-1\le n_0-2.
\]
By Claim~\ref{l:frontier} and the monotonicity of $\FL$, every admissible pair satisfies
\[
 1\le t\le s\le\FL(t)\le\FL(1)\le n_0-2.
\]
Consequently, the constants and order threshold in Lemma~\ref{l:expansion} can be chosen uniformly over all admissible pairs, and the constant implicit in its $\Oh(n^{-4})$ remainder depends only on $H$.

In the expansion \eqref{l:asymptotic-radius} supplied by Lemma~\ref{l:expansion}, the leading correction is proportional to $s^2+t^2$, with a positive coefficient of order $\Theta(n^{-2})$.
A strict difference in this integer quantity therefore dominates the remaining lower-order terms.

Among candidates with the same maximum sum of squares, the next correction is proportional to $st$, with a positive coefficient of order $\Theta(n^{-3})$.
A strict difference in $st$ dominates the $\Oh(n^{-4})$ error.
This proves the lexicographic rule.

If two positive pairs have the same sum of squares and the same product, they also have the same sum, since
$
 (s+t)^2=s^2+t^2+2st.
$
Their entries are therefore the roots of the same quadratic polynomial
$
 z^2-(s+t)z+st.
$
As both pairs are ordered with the larger entry first, they coincide.
Hence the maximizing pair is unique.

Finally, the coefficients used in these comparisons are positive for both odd and even $n$.
Thus the parameter rule is independent of the parity of $n$.
\end{proof}

The following criterion identifies exactly when the optimal nested parameter pair is diagonal.

\begin{corollary}\label{l:diagonal-test}
For sufficiently large $n$, the unique extremal graph is $\Ln kk$ if and only if
\begin{equation}\label{l:diagonal-condition}
 \FL(t)^2+t^2\le2k^2
 \qquad
 \text{for all }1\le t<k.
\end{equation}
A sufficient condition is the existence of a critical edge attaining $\beta_{u_1v_1}-1=k$ with $\alpha_{u_1v_1}=1$.
\end{corollary}

\begin{proof}
By Claim~\ref{l:frontier} in the proof of Theorem~\ref{l:parameters}, the diagonal pair $(k,k)$ is always a boundary candidate.
If another candidate has a strictly larger sum of squares, it has a strictly larger spectral radius for large $n$.

Suppose another candidate $(s,t)$ has sum of squares $2k^2$.
The arithmetic--geometric mean inequality gives
\[
 st\le\frac{s^2+t^2}{2}=k^2,
\]
with equality if and only if $s=t=k$.
Thus $(k,k)$ has a strictly larger product than every distinct candidate with the same sum of squares.
This proves the equivalence.

Under the stated sufficient condition, the definition gives $\FL(t)=k$ for every $t\ge1$.
Consequently, \eqref{l:diagonal-condition} holds.
\end{proof}

\begin{figure}[htbp]
\centering
\begin{tikzpicture}[
    x=0.75cm,
    y=0.75cm,
    line width=0.45pt,
    vertex/.style={circle,fill=black,inner sep=1.6pt}
]
\node[vertex] (la) at (-1.2,0) {};
\node[vertex] (lx) at (-3.2,0) {};
\node[vertex] (ra) at ( 1.2,0) {};
\node[vertex] (rx) at ( 3.2,0) {};

\foreach \i/\y in {1/-0.9,2/-0.3,3/0.3,4/0.9}{
    \node[vertex] (l\i) at (-2.2,\y) {};
    \node[vertex] (r\i) at ( 2.2,\y) {};
    \draw (la)--(l\i)--(lx);
    \draw (ra)--(r\i)--(rx);
}

\node[vertex] (u) at (0,1.2) {};
\draw (la)--(u)--(ra);
\draw (l1)--(r1);
\end{tikzpicture}
\caption{A counterexample for Corollary~\ref{l:diagonal-test}}
\label{fig:l-sufficient-not-necessary}
\end{figure}
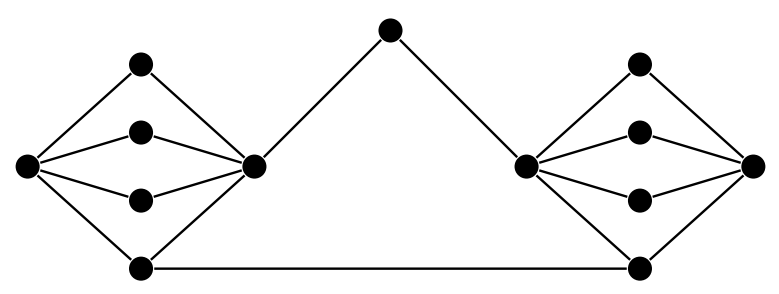

Figure \ref{fig:l-sufficient-not-necessary} (found through an AI-assisted search) shows that the sufficient condition in
Corollary~\ref{l:diagonal-test} is not necessary.

The next special case identifies when the extremal graph is a triangle attached to a balanced bipartite core.

\begin{corollary}\label{l:one-one}
For sufficiently large $n$, the unique extremal graph is $\Ln11$ if and only if $H$ has a critical edge whose endpoints both have degree two.
\end{corollary}

\begin{proof}
Such an edge gives $\beta_{u_1v_1}=2$, so $k=1$.
By Claim~\ref{l:frontier} in the proof of Theorem~\ref{l:parameters}, the only boundary pair is then $(1,1)$.

Conversely, if no critical edge has two degree-two endpoints, then $\beta_{u_1v_1}\ge3$ for every critical edge, and hence $k\ge2$.
The admissible pair $(k,k)$ has sum of squares at least $8$, strictly larger than that of $(1,1)$.
Thus $(1,1)$ cannot be optimal.
\end{proof}

\section{Proof of the main theorem, applications, and counterexamples}
\label{sec:applications}

This section completes the proof of the main theorem, illustrates the classification for all color-critical generalized theta graphs, and shows how asymmetric core-neighbor counts can separate spectral and edge extremality.

\begin{proof}[Proof of Theorem~\ref{thm:main}]
We first reduce the problem to connected forbidden graphs. Deleting a critical edge of $H$ makes every component
bipartite. Thus $H$ has a unique three-chromatic component $H_0$,
and we may write
\[
 H=H_0\cup H',
\]
where $H'$ is the union of the remaining components
and is therefore bipartite.

We claim that, for sufficiently large $n$,
\[
 \spex_3(n,H)=\spex_3(n,H_0),
 \qquad
 \SPEX_3(n,H)=\SPEX_3(n,H_0).
\]
If $H'$ is empty, this is immediate.
Otherwise, let $m_0=|V(H_0)|$.
Consider the graph
\[
 \Ttwo{n-2m_0-1}\cup C_{2m_0+1},
\]
where the two graphs have disjoint vertex sets and no
edges between them.
This graph is non-bipartite and $H_0$-free.
Indeed, the bipartite component cannot contain $H_0$,
whereas the cycle component has no odd cycle of length
at most $m_0$, and $H_0$ contains such an odd cycle.
Its spectral radius is at least $n/3$ for sufficiently
large $n$.

Every $H_0$-free graph is $H$-free.
Consequently, if $G\in\SPEX_3(n,H)$, then
$\rh(G)\ge n/3$.
The bound $\rh(G)^2\le2e(G)$ gives $e(G)\ge\frac{n^2}{18}.$
Suppose that $G$ contains a copy of $H_0$ on a vertex
set $W$.
Then
\[
 e(G-W)\ge\frac{n^2}{18}-m_0n.
\]
Choose a fixed integer $a\ge2$ such that
$H'\subseteq K_{a,a}$.
By the K\H{o}v\'ari--S\'os--Tur\'an
theorem~\cite{KovariSosTuran1954},
\[
 \operatorname{ex}(N,K_{a,a})
 =
 O_a(N^{2-1/a})
 =
 o(N^2).
\]
Hence $G-W$ contains $H'$ for sufficiently large $n$.
Together with the copy of $H_0$ on $W$, this gives a
copy of $H$ in $G$, a contradiction.
Therefore every graph in $\SPEX_3(n,H)$ is $H_0$-free.
Combining this conclusion with the inclusion of the
$H_0$-free admissible class in the $H$-free admissible
class proves the claimed equalities of the extremal
values and families.

The reduction also preserves all data used in the
classification.
Every critical edge lies in $H_0$, and its endpoint
degrees are unchanged.
Thus
\[
 \Estar(H)=\Estar(H_0),
 \qquad
 \dst(H)=\dst(H_0),
 \qquad
 \Pprof_L(H)=\Pprof_L(H_0).
\]
Moreover, $H$ contains a triangle if and only if $H_0$
does, because $H'$ is bipartite.

If $H-v$ is bipartite, then $v\in V(H_0)$ and all
neighbors of $v$ lie in $H_0$.
Every ordered bipartition of $H_0-v$ extends over the
bipartite components of $H'$, while every ordered
bipartition of $H-v$ restricts to one of $H_0-v$.
Consequently,
\[
 \Pprof_J(H)=\Pprof_J(H_0).
\]
The parameter rules in Subsections~\ref{sec:Jparameters} and~\ref{sec:Lparameters}
therefore give the same parameters for $H$ and $H_0$.
It suffices to prove the theorem for connected $H$.

Assume henceforth that $H$ is connected.
Let $G\in\SPEX_3(n,H)$ and let $\rho=\rh(G)$.
Lemma~\ref{prelim:lower} supplies an admissible graph
in each case and gives $\rho>(n-3)/2$.
Lemma~\ref{prelim:reduction} consequently implies that
$G$ is connected and that deleting a vertex $u$ of
minimum Perron coordinate leaves a bipartite graph.
We now distinguish the three cases in the statement.

\medskip
\noindent\emph{Case 1: $\dst(H)\ge3$.}
Theorem~\ref{j:complete-core} shows that $G-u=\Ttwo{n-1}$ and that both core neighborhoods of $u$ are nonempty and have sizes at most $m=\dst(H)-2$.
The orientation comparison in the proof of Theorem~\ref{j:exact-parameters} places the larger neighborhood in a smallest core part, so $G\cong\Jn st$ with $1\le t\le s\le m$.

By Lemma~\ref{j:frontier}, the graph $\Jn st$ is $H$-free if and only if $s\le\FJ(t)$.
If $s<\FJ(t)$, adding an allowed edge at the exceptional vertex gives a connected non-bipartite $H$-free graph with strictly larger spectral radius.
Thus an extremal graph must have $s=\FJ(t)$ and $t\in I_J$.
Conversely, every such boundary pair gives an admissible graph.
The comparison in Theorem~\ref{j:exact-parameters} first maximizes $\FJ(t)+t$ and then selects uniquely $t_J=\max T_J$ for odd $n$ and $t_J=\min T_J$ for even $n$.
Hence $G\cong\Jn{s_J}{t_J}$, where $s_J=\FJ(t_J)$.
Since an extremal graph exists and every extremal graph has this form, the displayed graph attains the maximum and is unique up to isomorphism.

\medskip
\noindent\emph{Case 2: $\dst(H)=2$ and $H$ contains a triangle.}
The graph $\ST{n-1}$ is non-bipartite and triangle-free, and is therefore $H$-free.
Theorem~\ref{s:structural} shows that every spectral extremal graph is triangle-free.
Lemma~\ref{prelim:triangle} then gives $\rho\le\rh(\ST{n-1}),$
with equality only when $G\cong\ST{n-1}$.
Admissibility gives the reverse inequality for the extremal value, proving both the asserted maximum and uniqueness.

\medskip
\noindent\emph{Case 3: $\dst(H)=2$ and $H$ is triangle-free.}
Theorem~\ref{l:family} gives $G\cong\Ln st$ with bounded integers $s\ge t\ge1$.
By the containment and boundary claims in the proof of Theorem~\ref{l:parameters}, admissibility is equivalent to $s\le\FL(t)$, and the possible values of $t$ are exactly $1,\ldots,k$.
As in Case 1, strict increase under the addition of an allowed edge forces $s=\FL(t)$ in an extremal graph.
Thus the optimization is over the finite list $(\FL(t),t)$, $1\le t\le k$.

Theorem~\ref{l:parameters} shows that these candidates are ordered first by $s^2+t^2$ and, among ties, by $st$.
The ordered maximizing pair is unique, because its sum of squares and product determine its two entries.
It is exactly $(s_L,t_L)$ defined in \eqref{l:selected-parameters}, independently of the parity of $n$.
Every boundary candidate is admissible, so this uniquely selected graph attains the maximum and
$\SPEX_3(n,H)=\{\Ln{s_L}{t_L}\}$.

Finally, Lemma~\ref{prelim:critical} gives $\dst(H)\ge2$, so the three cases exhaust all possibilities.
All reductions and comparisons require only finitely many order thresholds depending on the fixed graph $H$; take their maximum.
This proves the assertion for every sufficiently large $n$, and the disconnected reduction extends it from $H_0$ to $H$.
The formula for $\spex_3(n,H)$ follows from its definition.
\end{proof}

We next explain how Theorem~\ref{thm:main} recovers the known results for books and extends the classification for theta graphs.

\begin{remark}\label{rem:books}
For a book $B_t$ consisting of $t$ triangles sharing one edge, where $t\ge2$, we have
\[
 \dst(B_t)=t+1,
 \qquad
 \Pprof_J(B_t)=\{(1,t),(t,1)\}.
\]
Corollary~\ref{j:diagonal} gives the unique extremal graph $\Jn{t-1}{t-1}$ for sufficiently large $n$, recovering the theorem of Liu and Miao~\cite{LiuMiao2025}.
For $t=1$, the forbidden graph is a triangle, and the extremal graph is $\ST{n-1}$.
\end{remark}

A generalized theta graph $\theta(p_1,\ldots,p_k)$ consists of $k$ internally vertex-disjoint paths between two common endpoints, where $p_i$ denotes the number of edges of the $i$th path.
We assume the graph is simple, so at most one $p_i$ equals one.
It is non-bipartite precisely when both path-length parities occur.

The next corollary gives a complete classification for generalized color-critical theta graphs.

\begin{corollary}\label{cor:generalized-theta}
Let $H=\theta(p_1,\ldots,p_k)$ be a simple non-bipartite generalized theta graph. If $k=2$, then $H$ is an odd cycle, and for sufficiently large $n$,
\[
 \SPEX_3(n,H)=
 \begin{cases}
 \{\ST{n-1}\},&p_1+p_2=3,\\
 \{\Ln11\},&p_1+p_2\ge5.
 \end{cases}
\]

For $k\ge3$, $H$ has a critical edge if and only if one of the two path-length parities occurs exactly once.
If this condition holds and $n$ is sufficiently large, then
\[
 \SPEX_3(n,H)=
 \begin{cases}
 \{\Jn{k-2}{k-2}\},
 &\ell=1,\\[2pt]
 \{\ST{n-1}\},
 &\ell=2\text{ and another path has length }1,\\[2pt]
 \{\Ln{k-1}{k-1}\},
 &\ell=2\text{ and no path has length }1,\\[2pt]
 \{\Ln11\},
 &\ell\ge3.
 \end{cases}
\]
where $\ell$ is the length of that unique path.
\end{corollary}

\begin{proof}
Every cycle is the union of two distinct paths between the common endpoints.
Such a cycle is odd exactly when the two paths have opposite parities.
The graph is two-degenerate, so its non-bipartiteness implies $\chi(H)=3$.

For $k\ge3$, deleting an edge from one path destroys all odd cycles exactly when that path is the sole member of its parity class.
Thus the critical edges are precisely the edges of the unique minority-parity path.

If $\ell=1$, the critical edge joins the two common endpoints, both of degree $k$.
Deleting either endpoint gives profile $(1,k-1)$ up to orientation.
Deleting an internal vertex of any other path leaves an odd cycle using the critical edge and a remaining path.
Therefore
\[
 \Pprof_J(H)=\{(1,k-1),(k-1,1)\},
 \qquad
 \dst(H)=k.
\]
Corollary~\ref{j:diagonal} gives $\Jn{k-2}{k-2}$.

If $\ell=2$, each critical edge joins a degree-two internal vertex to a degree-$k$ endpoint.
A triangle occurs exactly when another path has length one.
In that case the middle line of the main theorem applies.
Otherwise every critical-edge profile is
\[
 (\alpha_{u_1v_1},\beta_{u_1v_1})=(1,k),
\]
so $\FL(t)=k-1$.
Corollary~\ref{l:diagonal-test} gives $\Ln{k-1}{k-1}$.

If $\ell\ge3$, the critical path contains an edge joining two degree-two internal vertices.
The graph is triangle-free, so Corollary~\ref{l:one-one} gives $\Ln11$.

When $k=2$, all edges of the odd cycle are critical, and all vertices have degree two.
The triangle and longer odd-cycle cases follow directly from the main theorem.
\end{proof}

Note that
\[
 \theta(1,2,2,2),\quad
 \theta(2,3,3,3),\quad
 \theta(3,2,2),\quad
 \theta(2,1,3)
\]
have extremal graphs
\[
 \Jn22,\qquad
 \Ln33,\qquad
 \Ln11,\qquad
 \ST{n-1},
\]
respectively.
The three-path specialization recovers the cases of Fang and Lin~\cite{FangLin2026}; the corollary also covers arbitrary numbers of paths whenever a critical edge exists.

The next example shows how parity resolves an edge-count tie between two admissible graphs in the $J$-family.


\begin{example}\label{ex:j-asymmetric}
Construct $H_J$ from two complete bipartite graphs with parts
\[
 \begin{aligned}
 V^{(1)}&=\{u_1^{(1)},u_3^{(1)}\},
 &V^{(2)}&=\{v_0^{(2)},v_1^{(2)},v_2^{(2)},v_3^{(2)}\},\\
 \widetilde V^{(1)}&=\{u_2^{(1)},u_4^{(1)},u_5^{(1)}\},
 &\widetilde V^{(2)}&=\{v_0^{(2)},v_4^{(2)},v_5^{(2)},v_6^{(2)}\},
 \end{aligned}
\]
sharing only $v_0^{(2)}$, and add the edge $u_1^{(1)}u_2^{(1)}$.
Then $H_J$ is connected, has minimum degree two, and has $u_1^{(1)}u_2^{(1)}$ as its only critical edge.
For sufficiently large $n$,
\[
 \SPEX_3(n,H_J)=
 \begin{cases}
 \{\Jn22\},&n\text{~is odd},\\
 \{\Jn31\},&n\text{~is even}.
 \end{cases}
\]
Both $\Jn22$ and $\Jn31$ are $H_J$-free and have the same number of edges, but their spectral radii are distinct.
Consequently,
\[
 \SPEX_3(n,H_J)\ne\operatorname{EX}_3(n,H_J).
\]

\smallskip\noindent
The graph has $12$ vertices and $21$ edges.
It contains the triangle $u_1^{(1)}u_2^{(1)}v_0^{(2)}$, and deleting $u_1^{(1)}u_2^{(1)}$ makes it bipartite.
Every other edge except $u_1^{(1)}v_0^{(2)},u_2^{(1)}v_0^{(2)}$ leaves this triangle intact.
Deleting $u_1^{(1)}v_0^{(2)}$ or $u_2^{(1)}v_0^{(2)}$ preserves, respectively, the odd cycle
$u_1^{(1)}v_1^{(2)}u_3^{(1)}v_0^{(2)}u_2^{(1)}u_1^{(1)}$
or
$u_2^{(1)}v_4^{(2)}u_4^{(1)}v_0^{(2)}u_1^{(1)}u_2^{(1)}$.
Thus $u_1^{(1)}u_2^{(1)}$ is the only critical edge, and $d_{H_J}(u_1^{(1)})=d_{H_J}(u_2^{(1)})=5$ gives $\dst(H_J)=5$.

Any vertex whose deletion makes $H_J$ bipartite must belong to $\{u_1^{(1)},u_2^{(1)},v_0^{(2)}\}$.
Deleting $u_1^{(1)}$ or $u_2^{(1)}$ gives profile $(1,4)$ up to orientation.
Deleting $v_0^{(2)}$ gives $(2,3)$ up to orientation, because $u_1^{(1)}u_2^{(1)}$ places $V^{(1)}$ and $\widetilde V^{(1)}$ in opposite bipartition classes.
Hence
\[
 \Pprof_J(H_J)=\{(1,4),(4,1),(2,3),(3,2)\},
 \qquad
 \FJ(1)=3,\quad \FJ(2)=2,\quad \FJ(3)=1.
\]
The feasible boundary pairs are $(3,1)$ and $(2,2)$, both of sum four.
Theorem~\ref{j:exact-parameters} gives the stated parity rule and strict spectral ordering.

Both graphs have $e(\Ttwo{n-1})+4$ edges.
If the edge and spectral extremal families were equal, the unique spectral extremal graph would be edge-extremal.
The other admissible graph, having the same edge count, would also be edge-extremal and hence spectral extremal, contradicting uniqueness.
\end{example}

The following example (found through an AI-assisted search) exhibits a strict separation between spectral and edge extremality.
In the $L$-family, edge count depends on $s+t$, whereas the leading spectral correction depends on $s^2+t^2$.
The following example makes these quantities select different pairs.

\begin{example}\label{ex:l-spectral-edge}
Let $H_L$ be obtained from the cycle $v_1v_2v_3v_4v_5v_1$ by attaching $0,3,1,1,3$ copies of $C_4$ at $v_1,\ldots,v_5$, respectively.
Each attachment identifies one vertex of the copy with the specified $v_i$; all other vertices are distinct.
Then $H_L$ is connected, triangle-free, has minimum degree two, and has no bridges. For sufficiently large $n$,
\[
 \SPEX_3(n,H_L)=\{\Ln72\}.
\]
However, $\Ln55$ is also $H_L$-free and
\[
 e(\Ln55)=e(\Ln72)+1.
\]
Consequently,
\[
 \SPEX_3(n,H_L)\cap\operatorname{EX}_3(n,H_L)
 =
 \varnothing.
\]

\smallskip\noindent
The graph has $29$ vertices and $37$ edges.
Its cycle vertices have degrees
\[
 (d(v_1),d(v_2),d(v_3),d(v_4),d(v_5))
 =
 (2,8,4,4,8),
\]
and every other vertex has degree two.
Every edge lies on a cycle, so there are no bridges.
The original $C_5$ is the only odd cycle, and each attached $C_4$ can be colored consistently with its shared vertex.
Thus $\chi(H_L)=3$.

The critical edges are exactly the five edges of the original cycle, with profiles
\[
 \begin{array}{c|c}
 \text{edge }v_iv_j&(\alpha_{v_iv_j},\beta_{v_iv_j})\\ \hline
 v_1v_2,\ v_5v_1&(1,8)\\
 v_2v_3,\ v_4v_5&(3,10)\\
 v_3v_4&(3,6)
 \end{array}
\]
and $\dst(H_L)=2$.
Therefore
\[
 k=5,
 \qquad
 \FL(1)=\FL(2)=7,
 \qquad
 \FL(t)=5\quad(3\le t\le5).
\]
The boundary candidates are
\[
 \begin{array}{c|c|c}
 (s,t)&s+t&s^2+t^2\\ \hline
 (7,1)&8&50\\
 (7,2)&9&53\\
 (5,3)&8&34\\
 (5,4)&9&41\\
 (5,5)&10&50
 \end{array}
\]
The sum of squares is uniquely maximized at $(7,2)$, so Theorem~\ref{l:parameters} gives $\SPEX_3(n,H_L)=\{\Ln72\}$.

Both $(7,2)$ and $(5,5)$ are admissible boundary pairs, and
\[
 e(\Ln72)=e(\Ttwo{n-2})+10,
 \qquad
 e(\Ln55)=e(\Ttwo{n-2})+11.
\]
Thus the unique spectral extremal graph has fewer edges than an admissible graph and cannot be edge-extremal.
For comparison, let $p=\lfloor(n-2)/2\rfloor$ and $q=\lceil(n-2)/2\rceil$.
Equation~\eqref{l:asymptotic-radius} gives
\[
 \rh(\Ln72)-\rh(\Ln55)
 =
 \frac{3q}{2(pq)^{3/2}}
 -
 \frac{11q}{(pq)^2}
 +
 O(n^{-4})
 >
 0
\]
for sufficiently large $n$.
\end{example}

\begin{figure}[htbp]
\centering
\begin{tikzpicture}[
  x=1cm,y=1cm,
  vertex/.style={
    circle,
    fill=black,
    inner sep=0pt,
    minimum size=3.2pt
  },
  every label/.style={font=\small},
  line width=0.55pt
]
  \node[vertex,label=above:$v_1$]
    (v1) at (0,1.6) {};
  \node[vertex,label=left:$v_2$]
    (v2) at (1.6,0.5) {};
  \node[vertex,label=above left:$v_3$]
    (v3) at (1,-1.3) {};
  \node[vertex,label=above right:$v_4$]
    (v4) at (-1,-1.3) {};
  \node[vertex,label=right:$v_5$]
    (v5) at (-1.6,0.5) {};

  \draw (v1)--(v2)--(v3)--(v4)--(v5)--(v1);

  \foreach \base/\ang/\idx in {
    v2/-45/1,v2/20/2,v2/85/3,
    v3/-65/4,v4/245/5,
    v5/95/6,v5/160/7,v5/225/8
  }{
    \begin{scope}[shift={(\base)},rotate=\ang]
      \node[vertex] (c\idx-1) at (0.65,-0.25) {};
      \node[vertex] (c\idx-2) at (1.30,0) {};
      \node[vertex] (c\idx-3) at (0.65,0.25) {};

      \draw
        (\base)--(c\idx-1)--(c\idx-2)--(c\idx-3)--(\base);
    \end{scope}
  }
\end{tikzpicture}
\caption{The graph $H_L$ in Example~4.4.}
\label{fig:HL}
\end{figure}
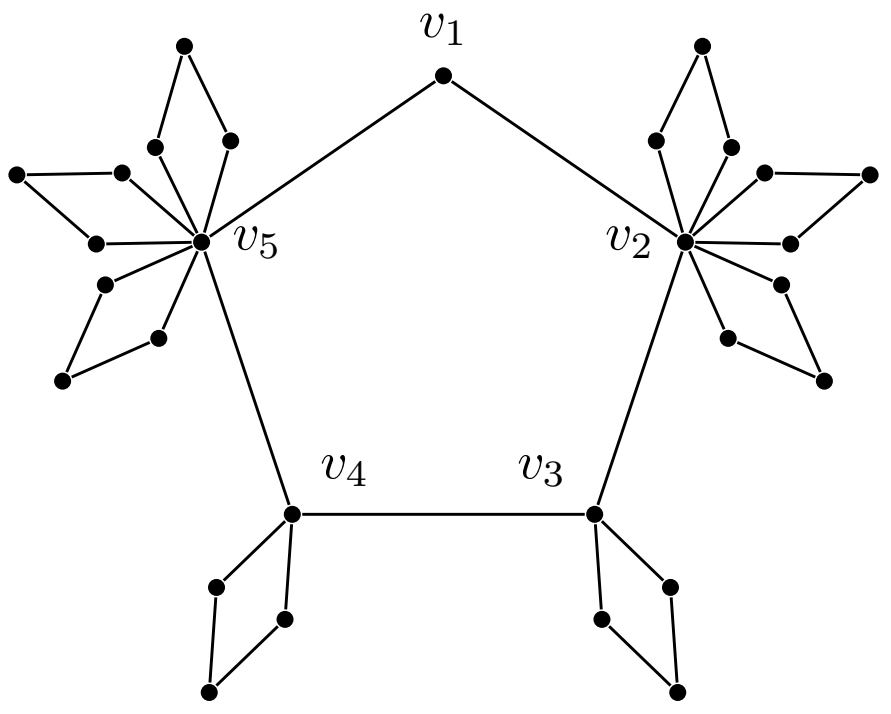

\begin{corollary}\label{rem:spectral-edge-conjecture}
The inclusion
$\SPEX_{r+1}(n,H)
 \subseteq\operatorname{EX}_{r+1}(n,H),$
conjectured in~\cite[Conjecture~1.1]{FangLin2026} for every fixed color-critical graph $H$ with $\chi(H)=r+1\ge3$ and sufficiently large $n$, is false.
It fails already for $r=2$ and a connected triangle-free forbidden graph with minimum degree two and no bridges.
\end{corollary}

\begin{proof}
The graph $H_L$ in Example~\ref{ex:l-spectral-edge} has all the stated properties, and its spectral and edge extremal families are disjoint.
\end{proof}


\begin{thebibliography}{99}

\bibitem{BhatiaDavis}
R.~Bhatia and C.~Davis,
A better bound on the variance,
\emph{Amer. Math. Monthly}
\textbf{107} (2000), 353--357.

\bibitem{BrouwerHaemers}
A.~E.~Brouwer and W.~H.~Haemers,
\emph{Spectra of Graphs},
Universitext, Springer, New York, 2012.

\bibitem{FangLin2026}
L.~Fang and H.~Lin,
Non-bipartite graphs without theta subgraphs,
\emph{J. Algebraic Combin.}
\textbf{63} (2026), Article~58.

\bibitem{HornJohnson}
R.~A.~Horn and C.~R.~Johnson,
\emph{Matrix Analysis},
second edition,
Cambridge University Press, Cambridge, 2013.

\bibitem{KovariSosTuran1954}
T.~K\H{o}v\'ari, V.~T.~S\'os, and P.~Tur\'an,
On a problem of K.~Zarankiewicz,
\emph{Colloq. Math.}
\textbf{3} (1954), 50--57.

\bibitem{LiLiuFeng2022survey}
Y.~Li, W.~Liu, and L.~Feng,
A survey on spectral conditions for some extremal
graph problems,
\emph{Adv. Math. (China)}
\textbf{51} (2022), 193--258.

\bibitem{LiPeng2023}
Y.~Li and Y.~Peng,
Refinement on spectral Tur\'an's theorem,
\emph{SIAM J. Discrete Math.}
\textbf{37} (2023), 2462--2485.

\bibitem{LinNingWu2021}
H.~Lin, B.~Ning, and B.~Wu,
Eigenvalues and triangles in graphs,
\emph{Combin. Probab. Comput.}
\textbf{30} (2021), 258--270.

\bibitem{LiuMiao2025}
R.~Liu and L.~Miao,
Spectral Tur\'an problem of non-bipartite graphs:
forbidden books,
\emph{European J. Combin.}
\textbf{126} (2025), Article~104136.

\bibitem{Nikiforov2007}
V.~Nikiforov,
Bounds on graph eigenvalues II,
\emph{Linear Algebra Appl.}
\textbf{427} (2007), 183--189.

\bibitem{Nikiforov2009}
V.~Nikiforov,
Spectral saturation: inverting the spectral Tur\'an theorem,
\emph{Electron. J. Combin.}
\textbf{16} (2009), Research Paper~33.

\bibitem{Nikiforov2011survey}
V.~Nikiforov,
Some new results in extremal graph theory,
in R.~Chapman (ed.),
\emph{Surveys in Combinatorics 2011},
London Mathematical Society Lecture Note Series,
vol.~392,
Cambridge University Press, Cambridge, 2011,
pp.~141--182.

\bibitem{OWu2026}
S.~O and J.~Wu,
A reduction principle for non-$r$-partite spectral
extremal problems, with a complete multipartite
classification,
\href{https://arxiv.org/abs/2607.00561v3}{arXiv:2607.00561v3}, 2026.

\bibitem{Simonovits1968}
M.~Simonovits,
A method for solving extremal problems in graph theory,
stability problems,
in \emph{Theory of Graphs
(Proc. Colloq., Tihany, 1966)},
Academic Press, New York, 1968, pp.~279--319.

\bibitem{Turan1941}
P.~Tur\'an,
On an extremal problem in graph theory,
\emph{Mat. Fiz. Lapok}
\textbf{48} (1941), 436--452 (in Hungarian).

\bibitem{WuLu2026}
Q.~Wu and Y.~Lu,
A counterexample to the Fang--Lin conjecture:
Edge and spectral extremality diverge near a Tur\'an graph,
\href{https://arxiv.org/abs/2609.37208v1}{arXiv:2609.37208v1}, 2026.

\bibitem{ZhaiLin2023}
M.~Zhai and H.~Lin,
A strengthening of the spectral chromatic
critical edge theorem: books and theta graphs,
\emph{J. Graph Theory}
\textbf{102} (2023), 502--520.

\bibitem{ZhangZhao2023}
Z.~Zhang and Y.~Zhao,
A spectral condition for the existence of cycles
with consecutive odd lengths in non-bipartite graphs,
\emph{Discrete Math.}
\textbf{346} (2023), Article~113365.

\bibitem{ZouFengLi2026}
L.~Zou, L.~Feng, and Y.~Li,
Spectral extremal problems for non-bipartite graphs
without odd cycles,
\emph{Discrete Math.}
\textbf{349} (2026), Article~114670.

\end{thebibliography}
\end{document}